\documentclass[11pt]{article}
\usepackage{smile}

\usepackage{fullpage}
\usepackage{framed}
\usepackage{mdframed}
\usepackage{enumerate}
\usepackage[inline]{enumitem}
\usepackage[T1]{fontenc}
\usepackage{moresize}
\usepackage{bm}
\usepackage{dsfont}
\usepackage{amsmath}
\usepackage{amssymb}
\usepackage{amsthm}
\usepackage{amsfonts}
\usepackage{stmaryrd}
\usepackage{array}
\usepackage{mathrsfs}
\usepackage{mathtools} 
\usepackage{extarrows}
\usepackage{stackrel}
\usepackage[nodisplayskipstretch]{setspace}
\usepackage{color}
\usepackage[usenames,dvipsnames]{xcolor}
\usepackage{cancel}
\usepackage{soul}
\usepackage{xfrac}
\usepackage{siunitx}
\usepackage{graphicx}
\usepackage{float}
\usepackage{rotating}
\usepackage{subcaption}
\usepackage{overpic}
\usepackage[all]{xy}
\usepackage{tikz}
\usetikzlibrary{arrows,matrix,positioning,calc,automata,patterns}
\usepackage{booktabs}
\usepackage{dcolumn}
\usepackage{multirow}
\usepackage{diagbox}
\usepackage{lscape}
\usepackage{rotating}
\usepackage{verbatim}
\usepackage{listings}
\usepackage{fancyvrb}
\usepackage{hyperref}
\usepackage[round]{natbib}
\usepackage{sectsty}

\usepackage[linesnumbered,ruled]{algorithm2e}
\usepackage{algpseudocode}
\SetKwInput{KwInput}{Input}             
\SetKwInput{KwInitialization}{Initialization}      
\SetKwInput{KWReturn}{Return}

\hypersetup{
    bookmarks=true,         
    unicode=false,          
    pdftoolbar=true,        
    pdfmenubar=true,        
    pdffitwindow=false,     
    pdfstartview={FitH},    
    pdftitle={My title},    
    pdfauthor={Author},     
    pdfsubject={Subject},   
    pdfcreator={Creator},   
    pdfproducer={Producer}, 
    pdfkeywords={key1, key2}, 
    pdfnewwindow=true,      
    colorlinks=true,        
    linkcolor=blue,         
    citecolor=blue,         
    filecolor=blue,         
    urlcolor=cyan           
}

\def\given{\,|\,}

\def\tr{\mathop{\text{tr}}\kern.2ex}

\def\P{{\mathrm P}}
\def\E{{\mathrm E}}
\def\R{{\mathbb R}}

\def\var{{\rm Var}}
\def\ind{{\mathds 1}}
\def\bP{{\mathrm P}}
\def\bE{{\mathrm E}}

\newcolumntype{L}[1]{>{\raggedright\let\newline\\\arraybackslash\hspace{0pt}}m{#1}}
\newcolumntype{C}[1]{>{  \centering\let\newline\\\arraybackslash\hspace{0pt}}m{#1}}
\newcolumntype{R}[1]{>{ \raggedleft\let\newline\\\arraybackslash\hspace{0pt}}m{#1}}
\newcolumntype{d}[1]{D{.}{.}{#1}}
\newcolumntype{H}{>{\setbox0=\hbox\bgroup}c<{\egroup}@{}}
\newcolumntype{Z}{>{\setbox0=\hbox\bgroup}c<{\egroup}@{\hspace*{-\tabcolsep}}}

\numberwithin{equation}{section}
\newtheorem{theorem}{Theorem}[section]
\newtheorem{lemma}{Lemma}[section]

\newtheorem{assumption}{Assumption}[section]

\newtheorem{definition}{Definition}[section]
\newtheorem{innercustomgeneric}{\customgenericname}

\providecommand{\customgenericname}{}
\newcommand{\newcustomtheorem}[2]{%
  \newenvironment{#1}[1]
  {%
   \renewcommand\customgenericname{#2}%
   \renewcommand\theinnercustomgeneric{##1}%
   \innercustomgeneric
  }
  {\endinnercustomgeneric}
}
\newcustomtheorem{customtheorem}{Theorem}
\newcustomtheorem{customassumption}{Assumption}
\newcustomtheorem{customlemma}{Lemma}
\newcustomtheorem{customexample}{Example}
\theoremstyle{definition}

\newtheorem{remark}{Remark}[section]

\newcommand{\dratio}{\frac{f_1(x)}{f_0(x)}}
\newcommand{\dif}{\mathrm{d}}
\newcommand{\bnorm}[1]{\Big\|#1\Big\|}

\let\hat\widehat
\let\bar\overline

\allowdisplaybreaks

\newcommand{\nb}[1]{{{#1}}}

\begin{document}

\setlength{\abovedisplayskip}{5pt}
\setlength{\belowdisplayskip}{5pt}
\setlength{\abovedisplayshortskip}{5pt}
\setlength{\belowdisplayshortskip}{5pt}
\hypersetup{colorlinks,breaklinks,urlcolor=blue,linkcolor=blue}

\title{\LARGE \nb{Limit theorems of matching estimators with a fixed number of matches}}

\author{
Songliang Chen\thanks{Department of Statistical Science, Duke University, Durham, NC 27708, USA; e-mail: \tt{songliang.chen@duke.edu}} ~~and~ Fang Han\thanks{Department of Statistics, University of Washington, Seattle, WA 98195, USA; e-mail: {\tt fanghan@uw.edu}}}

\date{\today}

\maketitle


\begin{abstract} 
\nb{This paper re-examines the limit theorems of Abadie and Imbens for nearest-neighbor matching estimators of average treatment effects with a fixed number of matches. We establish, for the first time, a non-normalized central limit theorem (CLT) with an explicitly calculated limiting variance. The key ingredients are to prove the convergence of the normalizing statistic appearing in the CLT of Abadie and Imbens to its mean, and to calculate the closed form of the limit of this mean. The former closes a gap in the argument of an unpublished work \citep{abadie2002simple}, while the latter resolves a question raised in \cite{abadie2006large}.}
\end{abstract}


\section{Introduction}

Consider estimating the population average treatment effect (ATE),
\[
\tau := \E\{Y(1)-Y(0)\},
\]
based on $n$ independent observations of $(Y(W),X,W)$, where $X\in\R^d$ represents some pre-treatment covariates, $W\in\{0,1\}$ indicates the treatment status, and $(Y(0),Y(1))\in\R^2$ are two potential outcomes \citep{neyman1923applications,rubin1974estimating} referring to being treated or not. 

For estimating $\tau$, the nearest neighbor (NN) matching estimator is intuitive and widely used \citep{stuart2010matching,imbens2024causal}. It imputes the missing potential outcomes by averaging the $M$ within-match outcomes in the opposite treatment group. In a landmark paper, Abadie and Imbens \citep[Theorem 4]{abadie2006large} proved the following central limit theorem for the {\it normalized} NN matching estimator, $\hat\tau_M$, with $M$ assumed to be fixed:
\[
V_{M}^{-1/2}\cdot \sqrt{n}(\hat\tau_M-B_{M}-\tau) \text{ converges in distribution to }\cN(0,1).
\]
Here $B_M=B_{M,n}$ represents the bias term and $V_M=V_{M,n}$ stands for a {\it random} approximation to the variance of $\sqrt{n}\hat\tau_M$. 

While techniques to correct the bias $B_{M}$ have been mature \citep{abadie2011bias}, still little is known about the random term $V_{M}$, whose stochastic behavior is ``difficult to work out'' \citep[Page 248 and Footnote 9]{abadie2006large}. \nb{Moreover, a non-normalized CLT for $\sqrt{n}(\hat\tau_M-B_{M}-\tau)$ itself has remained unavailable, making the theoretical justification of resampling procedures such as the $m$-out-of-$n$ bootstrap nontrivial.} The following theorem, presented informally here and to be rigorously stated in Section \ref{sec:theory}, settles \nb{both issues} and is our central result.

\begin{theorem}[Main theorem, informal]\label{thm:informal} Under certain regularity conditions and assuming $M$ to be fixed, the following hold true.
 \begin{enumerate}[label=(\roman*)]
\item The limit of $\E[V_{M}]$ exists and admits a closed form $\sigma_{M,d}^2$ to be introduced in Theorem \ref{thm4.2} ahead.
 \item \nb{$V_M$ converges in probability to $\E[V_M]$.}
\item $\sqrt{n}(\hat\tau_M-B_M-\tau)$ converges in distribution to $\cN(0,\sigma_{M,d}^2)$.
\end{enumerate}
\end{theorem}

\subsection{Related literature}

This paper contributes to a growing body of literature on large-sample theory for NN matching estimators. This line of research was pioneered by Abadie and Imbens in their series of seminal papers \citep{abadie2006large,abadie2008failure,abadie2011bias,abadie2012martingale,abadie2016matching} and has been extended by subsequent works, including \cite{lin2022regression}, \cite{lin2023estimation},  \cite{lu2023flexible}, \cite{huo2023adaptation}, \cite{cattaneo2023rosenbaum}, \cite{demirkaya2024optimal}, \cite{he2024propensity}, \cite{holzmann2024multivariate}, \cite{li2024matching}, \cite{ulloa2024propensity}, and \cite{lin2024consistency}, among many others. In particular, \cite{lin2023estimation} calculated the limiting variance of $\hat\tau_M$ as $M$ increases to infinity with $n$, and showed that the limit matches the corresponding semiparametric efficiency lower bound \citep{hahn1998role}. In this regard, the current paper extends their theory to the fixed-$M$ asymptotic regime.

For deriving the form of the limit of $\E[V_M]$ in Theorem \ref{thm:informal}(i), a crucial step is to calculate the limiting second moment of the volume of a stochastic object, which Abadie and Imbens termed the ``catchment area'' \citep[Page 260]{abadie2006large}. As pointed out in \citet[footnote 9]{abadie2006large} and \citet[Remark 2.1]{lin2023estimation}, this problem is closely related to a stochastic geometry problem of calculating the volume of Voronoi cells \citep{voronoi1908nouvelles}, and reduces to it when $M=1$. In this regard, our work also extends a result of Devroye, Gy\"{o}rfi, Lugosi, and Walk \citep[Theorem 1]{devroye2017measure} --- they calculated the limiting second moment of Voronoi's cells --- to the case when $M>1$. 

\nb{After deriving the limit of $\E[V_M]$, the remaining technical issue is to establish that $V_M$ converges in probability to this limit, as stated in Theorem \ref{thm:informal}(ii). An unpublished manuscript \citep{abadie2002simple} provided some heuristics toward this conclusion, but a fully rigorous argument has remained lacking. The present work fills this gap.}

More broadly speaking, this paper constitutes to the literature of statistical theory for NN-based methods \citep{MR3445317}, which aim to estimate a certain functional of the data generating distribution using a data-based NN graph. In this respect, the current paper is particularly related to a line of research on the asymptotic distribution-free properties of NN functionals \citep{MR532236,MR682809,shi2021ac,han2022azadkia}, to which this paper contributes a new distribution-free theorem (Theorem \ref{assump3.1} ahead, with two densities therein set to be identical).

\subsection{Notation and set-up}

For any positive integer $n$, denote $[n]=\{1,2,\ldots,n\}$. Let $\{(Y_i(0),Y_i(1),X_i,W_i)\}_{i\in[n]}$ be $n$ realizations of $(Y(0),Y(1),X,W)$. Write $Y=Y(W)$ and denote the support of $X$ by $\cX\subset \R^d$. Let $\{(Y_i=Y_i(W_i),X_i,W_i)\}_{i\in[n]}$ be the observation and  
\[
n_1:=\sum_{i=1}^nW_i~~~{\rm and}~~~n_0:=\sum_{i=1}^n(1-W_i)
\]
be the sizes of treated and control groups, respectively. Following \cite{abadie2006large}, we define the conditional average treatment effect
\[
\tau(x)=\E[Y(1)-Y(0) \given X=x]
\]
and introduce, for any $x\in\cX$ and $w\in \{0,1\}$, the following two conditional moments,
\[
\mu_w(x)=\E[Y\given X=x, W=w]~~~{\rm and}~~~\sigma_w^2(x)=\var(Y\given X=x,W=w).
\]
It is immediate that $\tau(x)=\mu_1(x)-\mu_0(x)$ if $(Y(0),Y(1))$ is independent of $W$ conditional on $X=x$ and $W$ is not degenerate given $X=x$. Furthermore, if $(Y(0),Y(1))$ is independent of $W$ conditional on almost all $X=x$ and $\P(W=1\given X)$ is bounded away from 0 and 1 almost surely, then 
\[
\tau=\E[\tau(X)]=\E[\mu_1(X)-\mu_0(X)].
\]
Lastly, for any $i\in[n]$, let $\varepsilon_i:=Y_i-\mu_{W_i}(X_i)$ be the residual and, for any $x\in\cX$, let $e(x):=\P(W=1\given X=x)$ be the propensity score \citep{rosenbaum1983central}.

\section{Nearest neighbor matching}

Following \cite{abadie2006large}, for any given number of matches $M$, we consider the following matching-based estimator of the ATE,
\[
\hat{\tau}_M = \frac{1}{n}\sum_{i = 1}^n\Big\{\hat{Y}_i(1) - \hat{Y}_i(0)\Big\},
\]
where for each $i\in[n]$, we define
\[
\hat{Y}_i(1)= \begin{cases}
        Y_i, & \mbox{ if } W_i=1,\\
        M^{-1}\sum_{j \in \cJ_M(i)}Y_j, & \mbox{ if } W_i=0
    \end{cases}
    ~~~{\rm and}~~~
\hat{Y}_i(0)= \begin{cases}
        Y_i, & \mbox{ if } W_i=0,\\
        M^{-1}\sum_{j \in \cJ_M(i)}Y_j, & \mbox{ if } W_i=1.
    \end{cases}
\]
Here $\cJ_M(i)$ indexes all $M$ NNs, based on the Euclidean metric $\norm{\cdot}$,  in the opposite treated group of subject $i$. In other words, $\cJ_M(i)$ indexes all $j \in [n]$ such that $W_j = 1 - W_i$ and 
\[
\sum_{k = 1, W_k = 1 - W_i}^n\ind\Big\{\bnorm{ X_k - X_i} \leq \bnorm{X_j - X_i} \Big\} \leq M.
\]

Let $K_M(i)$ denote the times the subject $i$ is used as a match, i.e.,
\[
K_M(i) := \sum_{j = 1, W_j = 1 - W_i}^n\ind\Big\{i \in  \cJ_M(j)\Big\}.
\]
The estimator $\hat{\tau}_M$ can also be written as 
\[
\hat{\tau}_M = \frac{1}{n}\sum_{i = 1}^n(2W_i -1)\Big\{1 + \frac{K_M(i)}{M}\Big\}Y_i,
\]
which, as correctly pointed out in \citet[Chapter 15.3.2]{ding2024first} and rigorously formalized in \cite{lin2023estimation}, resembles inverse probability weighted estimators.

\section{Theory}\label{sec:theory}

Our main theorem, Theorem \ref{thm:informal}, is proved under the following four sets of assumptions.

\begin{assumption} \label{assump4.1}
$\{(Y_i(0),Y_i(1),X_i,W_i)\}_{i\in[n]}$ are independently  drawn from $(Y(0),Y(1),X,W)$.
\end{assumption}

\begin{assumption} \label{assump4.2} For almost all $x\in\cX$, 
    \begin{itemize}
        \item [(i)] $(Y(0),Y(1))$ is independent of $W$ conditional on $X=x$;
        \item [(ii)] there exists some constant $\eta > 0$ such that $\eta < e(x) < 1 - \eta$.
    \end{itemize}
\end{assumption}

\begin{assumption} \label{assump4.3}
    \begin{itemize}
        \item [(i)] The support of $X$ is compact and convex;
        \item [(ii)]  the random vector $X$ has Lebesgue density $f_X$ and  there exist some constants $c_{\ell},  c_u\in (0,\infty)$ such that $c_{\ell} \leq \inf_{x\in\cX}f_X(x)\leq \sup_{x\in\cX}f_X(x) \leq c_u$.
    \end{itemize}
\end{assumption}

\begin{assumption}\label{assump4.4}
For $w = 0, 1$, 
 \begin{itemize}
     \item [(i)] $\mu_w(x)$ and $\sigma_w^2(x)$ are Lipschitz in $\cX$;
     \item [(ii)] $\E[Y^4 \given X = x, W = w]$ is uniformly bounded in $\cX$;
     \item [(iii)] There exist some constants $\overline{\sigma}^2, \underline{\sigma}^2 > 0$ such that $\underline{\sigma}^2 < \sigma_w^2(x) < \overline{\sigma}^2$ for all $x\in\cX$.    
\end{itemize}
\end{assumption}

The above four are identical to Assumptions 1-4 in \cite{abadie2006large}. 
In \cite{abadie2006large}, it is shown that $\hat{\tau}_M$ has the following decomposition
\[
\hat{\tau}_M - \tau = \overline{\tau(X)}  + E_M + B_M - \tau,
\]
with
\[
\overline{\tau(X)} := \frac{1}{n}\sum_{i = 1}^n\Big\{\mu_1(X_i)  - \mu_0(X_i)\Big\},~~
E_M :=  \frac{1}{n}\sum_{i = 1}^n(2W_i - 1)\Big\{1 + \frac{K_M(i)}{M}\Big\}\varepsilon_i,
\]
and
\[
B_M =  \frac{1}{n}\sum_{i = 1}^n(2W_i - 1)\Big[\frac{1}{M}\sum_{j \in \cJ_M(i)}\Big\{\mu_{1 - W_i}(X_i) - \mu_{1- W_i}(X_j)\Big\}\Big].
\]
Defining
\[
V^E =  \frac{1}{n}\sum_{i = 1}^n\Big\{1 + \frac{K_M(i)}{M}\Big\}^2\sigma_{W_i}^2(X_i),~~ 
V^{\tau(X)} = \E[(\tau(X) - \tau)^2],~~{\rm and}~~V_M:=V^E+V^{\tau(X)},
\]
we summarize the asymptotic properties of $\hat{\tau}_M$, derived in \cite{abadie2006large} and \cite{abadie2011bias}, as follows.

\begin{theorem}[\cite{abadie2006large, abadie2011bias}] \label{thm4.1} 
    \begin{itemize}
    \item[(i)]
Assume Assumption \ref{assump4.1}-\ref{assump4.4} hold. We then have
\[
   V_M^{-1/2}\cdot \sqrt{n}(\hat{\tau}_M - B_M - \tau) \text{ converges in distribution to } \cN(0, 1).
\]
\item[(ii)] Under additional smoothness conditions on $\mu_w$'s \citep[Theorem 2]{abadie2011bias}, there exists a bias-correction statistic $\hat B_M$ such that $\sqrt{n}(\hat B_M-B_M)$ converges in probability to 0 and 
\[
V_M^{-1/2}\cdot \sqrt{n}(\hat{\tau}_M^{\rm bc} - \tau) \text{ converges in distribution to } \cN(0, 1),
\]
with $\hat{\tau}_M^{\rm bc}:=\hat\tau_M-\hat B_M$.
\end{itemize}
\end{theorem}

It is evident that their theory focuses on the normalized $\hat\tau_M$, whose limiting variance relies on the asymptotic behavior of the statistic $V^E$. In \cite{abadie2006large}, the limit of $\E[V^E]$ is derived only in the special case of $d = 1$; see, for example, Theorem 5 therein. \nb{Moreover, the convergence of $V^E$ to $\E[V^E]$ was only briefly discussed in an unpublished manuscript \citep{abadie2002simple}, and its argument appears to be not fully rigorous.}

Using our \nb{variance analysis of $V^E$}, together with our analysis of high-order Voronoi cells \nb{to be detailed in Section \ref{voronoi} below}, we derive, for the first time, \nb{both the convergence of $V^E$ to $\E[V^E]$ and the limiting variance of $\hat\tau_M$} for a general dimension $d\geq 1$. The result is formalized in the following theorem, under an additional continuity assumption on the density functions. \nb{Note that our variance calculation agrees with the corresponding result in \citet[Theorem 5]{abadie2006large} when $d=1$. }

\begin{assumption}\label{assump4.5} For $w=0,1$, the Lebesgue density function of $X \given W=w$ is continuous.
\end{assumption}

\begin{theorem}[Main theorem] \label{thm4.2}   Assume Assumptions \ref{assump4.1}-\ref{assump4.5} hold.
\begin{itemize}
    \item[(i)]  We have 
\[
 \sqrt{n}(\hat{\tau}_M - B_M - \tau) \text{ converges in distribution to } \cN(0, \sigma^2_{M,d}),
 \]
 where $\sigma_{M,d}^2$ takes the form
 \begin{align*}
&\sigma^2_{M,d}:=\lim_{n\to\infty} \E[V_{M}]  \\
&=V^{\tau(X)} + \frac{1}{M^2}\E\Big[\sigma_1^2(X)\Big\{\frac{\alpha(M, d)}{e(X)} + \Big(\alpha(M, d) - M^2 - M\Big)e(X) + \Big(2M^2 + M - 2\alpha(M, d)\Big)\Big\}\Big] \\
    &+ \frac{1}{M^2}\E\Big[\sigma_0^2(X)\Big\{\frac{\alpha(M, d)}{1 - e(X)} + \Big(\alpha(M, d) - M^2 - M)(1 - e(X)\Big) + \Big(2M^2 + M - 2\alpha(M, d)\Big)\Big\}\Big],
\end{align*}
with $\alpha(M,d)$, defined in Equation \eqref{eq:alphaMd} ahead, being a {\it distribution-free} constant only depending on $M$ and $d$.
 \item[(ii)] Assume further the same conditions as in \citet[Theorem 2]{abadie2011bias}, the same bias-corrected estimator $\hat\tau_{M}^{bc}$ satisfies
 \[
  \sqrt{n}(\hat{\tau}_M^{\rm bc} - \tau) \text{ converges in distribution to } \cN(0, \sigma^2_{M,d}).
 \]
 \end{itemize}
\end{theorem}

\nb{
\begin{remark}\label{remark:conjecture}
We would like to highlight two interesting open questions. First, we conjecture that, for any $d\geq 1$, $\sigma^2_{M,d}$ converges to Hahn's semiparametric efficiency lower bound as $M\to\infty$. This is strongly suggested by \citet{lin2023estimation} and is known to hold when $d=1$ \citep[Theorem 5]{abadie2006large}. Second, for any fixed $M$, we conjecture that the term $\alpha(M,d)$ in Theorem \ref{thm4.2} is nonincreasing in $d$, as suggested by the numerical calculations in Tables \ref{tab:1} and \ref{tab:2} below. At present, we do not know how to prove either conjecture.
\end{remark}
}

\section{Measures of high-order Voronoi cells} \label{voronoi}

\nb{A key step in deriving the closed-form limiting variance in} Theorem \ref{thm4.2} is to quantify the volume of an NN-related stochastic object, which itself deserves a separate section to discuss. 

\nb{In the following, we consider a general situation as \citet[Section 2]{lin2023estimation}, and} let $\nu_0$ and $\nu_1$ be two laws on $\R^d$ with Lebesgue densities $f_0$ and $f_1$, respectively.  Let $X_1, \dots, X_{n}$ be $n$ independent draws from $\nu_0$.  

\begin{definition}[$M$-th NN map] Let $\cX_M(\cdot): \R^d \to \{X_i\}_{i = 1}^{n}$ be the map that returns the input $z$'s $M$-th NN in $\{X_i\}_{i = 1}^{n}$, that is, the value $x \in \{X_i\}_{i = 1}^{n}$ such that
\[
\sum_{i = 1}^{n}\ind\Big\{\bnorm{ X_i - z} \leq \bnorm{ x - z} \Big\} = M.
\]
\end{definition}

\begin{definition}[Catchment area] \label{def3.1} Let $\cA_M(\cdot): \bR^d \to \cB(\bR^d)$ be the  map from $\bR^d$ to the class of all Borel sets in $(\bR^d,\|\cdot\|)$ such that 
\[
\cA_M(x) = \cA_M\Big(x; \{X_i\}_{i = 1}^{n}\Big) := \Big\{z \in \bR^d: \bnorm{ z - x} \leq \bnorm{ \cX_M(z) - x} \Big\}. 
\]
\end{definition}\label{def3.2}

As noted in \citet[Remark 2.1]{lin2023estimation}, the sets $\{\cA_1(X_i)\}_{i\in[n]}$ are almost surely disjoint, partition $\R^d$ into $n$ polygons, and correspond exactly to the Voronoi cells introduced by \cite{voronoi1908nouvelles}. Accordingly, in this paper, we refer to the general sets $\{\cA_M(X_i)\}_{i\in[n]}$, for some $M\geq 1$, as {\it high-order Voronoi cells}. However, it should be noted that for $M>1$, $\{\cA_M(X_i)\}_{i\in[n]}$ are generally not disjoint. Additionally, they differ from the $K$-th order Voronoi cell/tessellation commonly defined in stochastic geometry; cf. \citet[Chapter 3.2]{boots1999spatial}.


Let $B(x,r)$ be the closed ball in $(\R,\|\cdot\|)$ of center $x\in\R^d$ and radius $r>0$. Following \cite{devroye2017measure}, introduce $D$ to be a random vector uniformly distributed in $B(0,1)$. Define $\bar{1} = (1, 0, \dots, 0) \in \bR^d$ and let $\bar{B} = B(\bar{1}, 1) \cup B(D, \Vert D \Vert)$. Define the random variable $T$ as 
\[
T = \frac{\lambda(\bar{B})}{\lambda(B(0,1))},
\]
with $\lambda$ representing the Lebesgue measure on $\bR^d$. Lastly, introduce
\[
\alpha(d) := \E\Big[\frac{2}{T^2}\Big].
\]

\nb{\citet[Theorem 1]{devroye2017measure} focuses on the case of $\nu_1(\cA_1(x))$ with $\nu_0=\nu_1$. Theorem \ref{thm3.1} below generalizes their result to the measures of high-order Voronoi cells, while further allowing $\nu_0$ and $\nu_1$ to be possibly distinct.}

\begin{assumption}\label{assump3.1} Assume $\nu_1$ is absolutely continuous with regard to $\nu_0$ so that the Lebesgue densities $f_0$ and $f_1$ share a common support $\cX$ that is compact. Furthermore, assume $f_0$ and  $f_1$ to be continuous in $\cX$.
\end{assumption}

\begin{theorem} \label{thm3.1} 
\begin{itemize}
\item[(i)] Under Assumption \ref{assump3.1}, for $\nu_0$-almost all $x$ we have
    \begin{align}
        \lim_{n\to\infty} n\E\Big[\nu_1(\cA_M(x)) \mid X_1 = x\Big] &= M\cdot \frac{f_1(x)}{f_0(x)}, \label{eq:thm3.1-1}\\
        \lim_{n\to\infty}n^2\E\Big[\nu_1(\cA_M(x))^2 \mid X_1 = x\Big] &= \alpha(M,d)\cdot \left(\frac{f_1(x)}{f_0(x)}\right)^2, \label{eq:thm3.1-2}
    \end{align}
    where we define
    \begin{align}\label{eq:alphaMd}
    \alpha(M, d) := \alpha(d)\sum_{i + k \leq M - 1, j +k \leq M - 1}\frac{c_{ijk}(d)(i + j + k + 1)!}{i!j!k!},
    \end{align}
    with $c_{ijk}(d)$ introduced in Lemma \ref{lem5.2} ahead being distribution-free constant depending only on $i, j, k,$ and $d$.
    \item[(ii)] For $d = 1$, we have $\alpha(M, 1) = M(2M+1)/2$. For $M=1$, we have $\alpha(1,d)=\alpha(d)$.
\end{itemize}    
\end{theorem}

In the special case where $\nu_0=\nu_1$, Theorem \ref{thm3.1} demonstrates that both $ n\E[\nu_0(\cA_M(x)) \mid X_1 = x]$ and $n^2\E[\nu_0(\cA_M(x))^2 \mid X_1 = x]$ have distribution-free limits that only depend on the NN number, $M$, and the dimension, $d$. While the closed form of $\alpha(M,d)$ is generally unavailable unless either $d$ or $M$ is 1, these values can be computed numerically. 

The results are presented in Tables \ref{tab:1} and \ref{tab:2}.

\begin{table}[H]\label{tab:1}
        \renewcommand{\arraystretch}{1.5}
                \caption{Values of $\alpha(d)$}
        \centering
        \nb{
        \begin{tabular}{cccccccccccc}
            \hline
            & $d= 1$ & $d = 2$ & $d = 3$ &$d = 4$ & $d = 5$ & $d = 6$ & $d = 7$ & $d = 8$ & $d = 9$ & $d = 10$ &\\
            \hline
             $\alpha(d)$& 1.50& 1.16  & 1.14 & 1.11& 1.08  & 1.06& 1.04 & 1.03 & 1.02 & 1.02 &\\
            \hline
        \end{tabular}
        }
        \label{tab:1}
    \end{table}

    \begin{table}[H]\
        \renewcommand{\arraystretch}{1.5}
                \caption{Values of $\alpha(M,d)$}
        \centering
        \nb{
        \begin{tabular}{ccccccccccc}
            \hline
             $\alpha(M, d) $& $d = 2$ & $d = 3$ &$d = 4$ & $d = 5$ & $d = 6$ & $d = 7$ & $d = 8$ & $d = 9$ & $d = 10$ &\\
            \hline
             $M = 2$& 4.34 & 4.31 & 4.24 & 4.17 & 4.12 & 4.08 & 4.05 & 4.03  & 4.01 & \\
             \hline
             $M = 3$& 9.55 & 9.50 & 9.38 & 9.27 & 9.19 & 9.13 & 9.07 & 9.03 &  9.00&\\
             \hline
              $M = 4$& 16.76 & 16.70 & 16.53 & 16.38 & 16.27 & 16.17 & 16.09 & 16.03 & 15.97 & \\
             \hline
              $M = 5$& 25.99 & 25.91 & 25.69 & 25.51 & 25.35 & 25.22 & 25.11 & 25.01 & 24.93&\\
             \hline
              $M = 6$&  37.22 & 37.13 & 36.87 & 36.63 & 36.44 & 36.27 & 36.13 & 35.99 & 35.87 &\\
             \hline
              $M = 7$&  50.45 & 50.36  & 50.05 & 49.77 & 49.53 & 49.32 & 49.14 & 48.96 & 48.81 &\\
             \hline
              $M = 8$&  65.69 & 65.60 & 65.24  & 64.92 & 64.63 & 64.38 & 64.15 & 63.92 &  63.72&\\
             \hline
              $M = 9$& 82.94 & 82.85 & 82.44 & 82.07 & 81.74 & 81.44 & 81.16  & 80.88 & 80.63 &\\
             \hline
              $M = 10$& 102.18 & 102.11 & 101.66 & 101.24 & 100.85 & 100.50 & 100.17 & 99.83 & 99.53 &\\
             \hline
        \end{tabular}
        }
        \label{tab:2}
    \end{table}

\section{Proofs of the main results}

\subsection{Proof of Theorem \ref{thm3.1}}

Let $Z_1, Z_2$ be two independent copies of $Z \sim \nu_1$. Define 
\begin{align*}
&V_1 = \nu_0(B(Z_1, \Vert Z_1 - x \Vert)),~~ V_2 = \nu_0(B(Z_2,\Vert Z_2 - x \Vert)),\\
~~ {\rm and}~~ &V = \nu_0\big(B(Z_1, \Vert Z_1 - x \Vert) \cup B(Z_2, \Vert Z_2 - x \Vert)\big).
\end{align*}
We first introduce the following lemma.
\begin{lemma}\label{lem5.1}
Denote the Lebesgue density functions of $V_1, V$ by $f_{V_1}, f_V$. We then have
\[
\begin{aligned}
\lim_{t \to 0+} f_{V_1}(t) &= \frac{f_1(x)}{f_0(x)} ~~~{\rm and}~~~\lim_{t \to 0+} \frac{f_V(t)}{t} &= \alpha(d)\Big\{\frac{f_1(x)}{f_0(x)}\Big\}^2.
\end{aligned}
\]
\end{lemma}

\noindent\textbf{Step 1.} We first prove \eqref{eq:thm3.1-1}. Observe that
\[
\begin{aligned}
&\quad \bE\Big[\nu_1(\cA_M(x)) \mid X_1 = x\Big]\\
&= \bP\Big(Z_1 \in \cA_M(x) \mid X_1 = x\Big)\\
&= \bP\Big(\Vert x - Z_1\Vert \leq \Vert \cX_M(Z_1) - Z_1 \Vert \mid X_1 = x\Big) \\
& = \bP\Big(\text{at most } M - 1 \text{ indices } i \in \{2, \dots, n\} \text{ satisfy } \Vert x - Z_1\Vert > \Vert X_i - Z_1\Vert  \mid X_1 = x\Big)\\
& = \sum_{0 \leq k \leq M - 1}\bP\Big(\text{there exist exactly } k \text{ indices } i \in \{2, \dots n\} \text{ such that } \Vert x - Z_1\Vert > \Vert X_i - Z_1\Vert \mid X_1 = x\Big)\\
& = \sum_{0 \leq k \leq M - 1} \binom{n - 1}{k} \bP\Big(\Vert x - Z_1\Vert > \Vert X_2 - Z_1 \Vert\Big)^k \cdot \bP\Big(\Vert x - Z_1\Vert \leq \Vert X_2 - Z_1 \Vert\Big)^{n - 1 - k}\\
& = \sum_{0 \leq k \leq M - 1} \binom{n - 1}{k} \bP\Big(X_2 \in B_{Z_1, \Vert Z_1 - x\Vert}\Big)^k \bP\Big(X_2 \notin B_{Z_1, \Vert Z_1 - x\Vert}\Big)^{n - k - 1}\\
& = \sum_{0 \leq k \leq M - 1} \bE\Big[\binom{n - 1}{k} V_1^k(1 - V_1)^{n - 1 - k}\Big].
\end{aligned}
\]
For each $0 \leq k \leq M - 1$ and all $0 < \delta < 1$, we have
\[
n\bE\Big[\binom{n - 1}{k}V_1^k(1 - V_1)^{n - 1 - k}\ind\{V_1 \geq \delta\}\Big] \leq n\binom{n - 1}{k}(1 - \delta)^{n - 1 - k}.
\]
Since $n\binom{n - 1}{k}$ is a polynomial of $n$ and $(1 - \delta)^n$ decays to 0 exponentially fast, for any fixed $0 < \delta < 1$, we then have
\[
n\bE\Big[\binom{n - 1}{k}V_1^k(1 - V_1)^{n - 1 - k}\ind\{V_1 \geq \delta\}\Big] \to 0,  \quad \text{ as } n \to \infty.
\]
By Lemma \ref{lem5.1}, for any $\varepsilon > 0$, there exists $\delta\in(0,1)$ such that for any $v_1 \in [0, \delta]$, the density $f_{V_1}(v_1)$ satisfies
\[
(1 - \varepsilon)\dratio \leq f_{V_1}(v_1) \leq (1 + \varepsilon)\dratio.
\]
Thus
\begin{align*}
    &\quad n\bE\Big[\binom{n - 1}{k}V_1^k(1 - V_1)^{n - 1 - k}\ind\{V_1 < \delta\}\Big] \\
    &= n\binom{n - 1}{k}\int_0^{\delta}v_1^k(1 - v_1)^{n - 1 - k}f_{V_1}(v_1) \dif v_1 \\
    & \leq (1 + \varepsilon)\dratio n\binom{n - 1}{k} \int_0^\delta v_1^k(1 - v_1)^{n - 1- k} \dif v_1 \\
    & =  (1 + \varepsilon)\dratio n\binom{n - 1}{k} \left(\int_0^1  v_1^k(1 - v_1)^{n - 1- k} \dif v_1 - \int_\delta^1 v_1^k(1 - v_1)^{n - 1- k} \dif v_1 \right) \\
    & = (1 + \varepsilon)\dratio n\binom{n - 1}{k} \frac{\Gamma(k + 1)\Gamma(n - k)}{\Gamma(n + 1)} - (1 + \varepsilon)\dratio n\binom{n - 1}{k}\int_\delta^1 v_1^k(1 - v_1)^{n - 1- k} \dif v_1 \\
    & = (1 + \varepsilon)\dratio - (1 + \varepsilon)\dratio n\binom{n - 1}{k}\int_\delta^1 v_1^k(1 - v_1)^{n - 1- k} \dif v_1 .
\end{align*}
Since 
\[
 n\binom{n - 1}{k}\int_\delta^1 v_1^k(1 - v_1)^{n - 1- k} \dif v_1  \leq n\binom{n - 1}{k}(1 - \delta)^{n - 1 - k} \to 0, \quad \text{as } n \to \infty,
\]
we obtain
\begin{align}\label{eq:5.1}
\limsup_{n \to \infty} n\bE\Big[\binom{n - 1}{k}V_1^k(1 - V_1)^{n - 1 - k}\Big] \leq (1 + \varepsilon)\dratio. 
\end{align}
On the other hand, we can similarly derive the lower bound since
\begin{align*}
    &\quad n\bE\Big[\binom{n - 1}{k}V_1^k(1 - V_1)^{n - 1 - k}\ind\{V_1 < \delta\}\Big] \\    
    & \geq (1 - \varepsilon)\dratio n\binom{n - 1}{k} \int_0^\delta v_1^k(1 - v_1)^{n - 1- k} \dif v_1 \\
    & =  (1 - \varepsilon)\dratio n\binom{n - 1}{k} \left(\int_0^1  v_1^k(1 - v_1)^{n - 1- k} \dif v_1 - \int_\delta^1 v_1^k(1 - v_1)^{n - 1- k} \dif v_1 \right) \\
    & = (1 - \varepsilon)\dratio n\binom{n - 1}{k} \frac{\Gamma(k + 1)\Gamma(n - k)}{\Gamma(n + 1)} - (1 + \varepsilon)\dratio n\binom{n - 1}{k}\int_\delta^1 v_1^k(1 - v_1)^{n - 1- k} \dif v_1 \\
     & = (1 - \varepsilon)\dratio - (1 + \varepsilon)\dratio n\binom{n - 1}{k}\int_\delta^1 v_1^k(1 - v_1)^{n - 1- k} \dif v_1 .
\end{align*}
Thus
\begin{align}\label{eq:5.2}
\liminf_{n \to \infty} n\bE\Big[\binom{n - 1}{k}V_1^k(1 - V_1)^{n - 1 - k}\Big] \geq (1 - \varepsilon)\dratio. 
\end{align}
Since $\varepsilon$ is taken arbitrarily, combining \eqref{eq:5.1} and \eqref{eq:5.2}, we obtain
\[
\lim_{n \to \infty} n\bE\Big[\binom{n - 1}{k}V_1^k(1 - V_1)^{n - 1 - k}\Big] = \dratio.
\]
Therefore we have
\[
\lim_{n \to \infty}n\bE\Big[\nu_1(\cA_M(x))\mid X_1 = x\Big] = \lim_{n \to \infty} \sum_{0 \leq k \leq M-1}n\bE\Big[\binom{n - 1}{k}V_1^k(1 - V_1)^{n - 1 - k}\Big] = M\dratio,
\]
which completes the proof of \eqref{eq:thm3.1-1}.\\

\noindent\textbf{Step 2.} We next prove \eqref{eq:thm3.1-2}. Similar to the previous argument, we observe
\begin{align*}
&\quad \bE\Big[\nu_1(\cA_M(x))^2 \mid X_1 = x\Big]\\
&= \bP\Big(Z_1 \in \cA_M(x), Z_2 \in \cA_M(x) \mid X_1 = x\Big)\\
&= \bP\Big(\Vert x - Z_1\Vert \leq \Vert \cX_M(Z_1) - Z_1 \Vert, \Vert x - Z_2\Vert \leq \Vert \cX_M(Z_2) - Z_2 \Vert  \mid X_1 = x\Big) .
\end{align*}
Define the following four random sets
\begin{align*}
    & A = \{l \in \{2, \dots, n\}: \Vert x - Z_1\Vert \leq \Vert X_l - Z_1\Vert, \Vert x - Z_2\Vert \leq \Vert X_l - Z_2\Vert \}, \\
    & B = \{l \in \{2, \dots, n\}: \Vert x - Z_1\Vert \leq \Vert X_l - Z_1\Vert, \Vert x - Z_2\Vert > \Vert X_l - Z_2\Vert \}, \\
    & C = \{l \in \{2, \dots, n\}: \Vert x - Z_1\Vert > \Vert X_l - Z_1\Vert, \Vert x - Z_2\Vert \leq \Vert X_l - Z_2\Vert \}, \\
    & D = \{l \in \{2, \dots, n\}: \Vert x - Z_1\Vert > \Vert X_l - Z_1\Vert, \Vert x - Z_2\Vert > \Vert X_l - Z_2\Vert \}.
\end{align*}
Given $Z_1, Z_2$, let $p_A, p_B, p_C, p_D$ be the probabilities of a given index $l$ belonging to $A, B, C, D$, respectively. It then holds true that
\begin{align*}
     p_A &:= \bP\Big(\Vert x - Z_1\Vert \leq \Vert X_l - Z_1\Vert, \Vert x - Z_2\Vert \leq \Vert X_l - Z_2\Vert \mid Z_1, Z_2\Big) \\
    & = \nu_0 \Big(B(Z_1, \Vert Z_1 - x\Vert)^c \cap B(Z_2, \Vert Z_2 - x\Vert)^c\Big)\\
    & = \nu_0\Big((B(Z_1, \Vert Z_1 - x\Vert) \cup B(Z_2, \Vert Z_2 - x\Vert))^c\Big)\\
    & = 1 - V,\\
    p_B &:= \bP\Big(\Vert x - Z_1\Vert \leq \Vert X_l - Z_1\Vert, \Vert x - Z_2\Vert > \Vert X_l - Z_2\Vert \mid Z_1, Z_2\Big) \\
    & = \nu_0 \Big(B(Z_1, \Vert Z_1 - x\Vert)^c \cap B(Z_2, \Vert Z_2 - x\Vert)\Big) \\
    & = V - V_1,\\
    p_C &:= \bP\Big(\Vert x - Z_1\Vert > \Vert X_l - Z_1\Vert, \Vert x - Z_2\Vert \leq \Vert X_l - Z_2\Vert \mid Z_1, Z_2\Big) \\
    & = \nu_0 \Big(B(Z_1, \Vert Z_1 - x\Vert) \cap B(Z_2, \Vert Z_2 - x\Vert)^c\Big) \\
    & = V - V_2,\\
    p_D & := \bP\Big(\Vert x - Z_1\Vert > \Vert X_l - Z_1\Vert, \Vert x - Z_2\Vert > \Vert X_l - Z_2\Vert \mid Z_1, Z_2\Big) \\
    & = \nu_0 \Big(B(Z_1, \Vert Z_1 - x\Vert)^c \cap B(Z_2, \Vert Z_2 - x\Vert)^c\Big)\\
    & = V_1 + V_2 - V.
\end{align*}
Observe that $A, B, C, D$ are pairwisely disjoint, $|A| + |B| + |C| + |D| = n - 1$, and we have
\begin{align*}
B\cup D &= \Big\{l \in \{2, \dots, n\}: \Vert x - Z_1\Vert > \Vert X_l - Z_1\Vert \Big\},\\
 C \cup D &= \Big\{l \in \{2, \dots, n\}: \Vert x - Z_2\Vert > \Vert X_l - Z_2\Vert \Big\}.
\end{align*}
We then have
\[
\begin{aligned}
    & \quad \bP\Big(\Vert x - Z_1\Vert \leq \Vert \cX_M(Z_1) - Z_1 \Vert, \Vert x - Z_2\Vert \leq \Vert \cX_M(Z_2) - Z_2 \Vert  \mid X_1 = x\Big) \\
    & = \bP\Big(|B \cup D| \leq M - 1, |C \cup D| \leq M - 1 \mid X_1 = x\Big) \\
    & = \bP\Big(|B| + |D| \leq M - 1, |C| + |D| \leq M - 1\Big) \\
    & = \sum_{i + k , j + k \leq M - 1}\bE\Big[\bP(|D| = k, |B| = i, |C| = j \mid Z_1, Z_2)\Big] \\
    & = \sum_{i + k , j + k \leq M - 1}\binom{n - 1}{k}\binom{n - 1 - k}{i}\binom{n - 1 - k - i}{j}\bE\Big[p_A^{n - 1 - i - j - k}p_B^ip_C^jp_D^k\Big] \\
    & =  \sum_{i + k , j + k \leq M - 1} \frac{(n - 1)!}{i!j!k!(n - 1 - i - j - k)!} \bE\Big[(V - V_1)^i(V - V_2)^j(V_1 + V_2 - V)^k(1 - V)^{n - 1 - i - j - k}\Big].
\end{aligned}
\]
Define 
\[
\begin{aligned}
&P_{ijk}(v, v_1, v_2) = \Big(1 - \frac{v_1}{v}\Big)^i\Big(1 - \frac{v_2}{v}\Big)^j\Big(\frac{v_1 + v_2}{v} - 1\Big)^k,\\
&Q_{ijk}(v) = \bE\Big[\Big(1 - \frac{V_1}{V}\Big)^i\Big(1 - \frac{V_2}{V}\Big)^j\Big(\frac{V_1 + V_2}{V} - 1\Big)^k\mid V = v\Big] = \bE\Big[P_{ijk}(V, V_1, V_2) \mid V = v\Big]  .
\end{aligned}
\]
We then have
\begin{align*}
    & \quad \bE\Big[(V - V_1)^i(V - V_2)^j(V_1 + V_2 - V)^k(1 - V)^{n - 1 - i - j - k}\Big] \\
    & = \bE\Big[\bE[(V - V_1)^i(V - V_2)^j(V_1 + V_2 - V)^k \mid V](1 - V)^{n - 1 - i - j - k}\Big] \\
    & = \bE\Big[Q_{ijk}(V)V^{i + j + k}(1 - V)^{n - 1 - i - j - k}\Big].
    \end{align*}
    

\begin{lemma}[Distribution-free limits] \label{lem5.2} It holds true that
\[
\lim_{v \to 0+}Q_{ijk}(v) = c_{ijk}(d),
\]
where $c_{ijk}(d)$ is a positive constant depending only on the indices $i, j, k$ and dimension $d$. Specifically, when the dimension $d = 1 $, this constant can be explicitly calculated as
\[
c_{ijk}(1) = \begin{cases}
        1, & \mbox{ if } i = j = k = 0,\\
       \frac{1}{3}\frac{i!j!}{(i + j + 1)!}, & \mbox{ if } i, j \neq 0, k = 0,\\
        \frac{1}{3}\frac{i!k!}{(i + k + 1)!}, & \mbox{ if } i, k \neq 0, j = 0,\\
        \frac{1}{3}\frac{j!k!}{(j + k + 1)!}, & \mbox{ if } j, k \neq 0, i = 0,\\
        \frac{2}{3}\frac{1}{k + 1}, & \mbox{ if } k \neq 0, i = j = 0,\\
        \frac{2}{3}\frac{1}{j + 1}, & \mbox{ if } j \neq 0, i = k = 0,\\
        \frac{2}{3}\frac{1}{i + 1}, & \mbox{ if } i \neq 0, j = k = 0,\\
       0, & \mbox{ if } i, j, k \neq 0.
    \end{cases}
\]

\end{lemma}

Back to the proof, Lemmas \ref{lem5.1} and \ref{lem5.2} combined imply that, for any $\varepsilon > 0$, there exists some $\delta > 0$ such that for any $0 < v < \delta$, we have
\[
(1 - \varepsilon)\alpha(d)\Big\{\dratio\Big\}^2 \leq f_V(v) \leq (1 + \varepsilon)\alpha(d)\Big\{\dratio\Big\}^2,
\]
\[
(1 - \varepsilon)c_{ijk}(d) \leq Q_{ijk}(v) \leq (1 + \varepsilon)c_{ijk}(d).
\]
Similar to the argument in {\bf Step 1}, we now have
\begin{align*}
    & \quad n^2\frac{(n - 1)!}{i!j!k!(n - 1 - i - j - k)!} \bE\Big[(V - V_1)^i(V - V_2)^j(V_1 + V_2 - V)^k(1 - V)^{n - 1 - i - j - k} \ind \{V \leq \delta\}\Big] \\
    & = n^2\frac{(n - 1)!}{i!j!k!(n - 1 - i - j - k)!} \bE\Big[Q_{ijk}(V)V^{i + j + k}(1 - V)^{n - 1 - i - j - k} \ind \{V \leq \delta \}\Big]\\
    & = n^2\frac{(n - 1)!}{i!j!k!(n - 1 - i - j - k)!} \int_0^\delta Q_{ijk}(v)v^{i + j + k}(1 - v)^{n - 1 - i - j - k} \dif v \\
    & \leq (1 + \varepsilon)^2  c_{ijk}(d) \alpha(d) \Big(\dratio\Big)^2\frac{n^2(n - 1)!}{i!j!k!(n - 1 - i - j - k)!} \int_0^\delta v^{i + j + k + 1}(1 - v)^{n - 1 - i - j - k} \dif v \\
    & = (1 + \varepsilon)^2  c_{ijk}(d) \alpha(d) \Big(\dratio\Big)^2\frac{n^2(n - 1)!}{i!j!k!(n - 1 - i - j - k)!} (\int_0^1 - \int_\delta^1) v^{i + j + k + 1}(1 - v)^{n - 1 - i - j - k} \dif v \\
    & = (1 + \varepsilon)^2  c_{ijk}(d) \alpha(d) \Big(\dratio\Big)^2\frac{n^2(n - 1)!}{i!j!k!(n - 1 - i - j - k)!} \frac{\Gamma(i + j + k + 2)\Gamma(n - i - j - k)}{\Gamma(n + 2)} \\
    & \quad - (1 + \varepsilon)^2  c_{ijk}(d) \alpha(d) \Big(\dratio\Big)^2\frac{n^2(n - 1)!}{i!j!k!(n - 1 - i - j - k)!} \int_\delta^1 v^{i + j + k + 1}(1 - v)^{n - 1 - i - j - k} \dif  v\\ 
    & = (1 + \varepsilon)^2   \alpha(d) \frac{n}{n + 1}\frac{c_{ijk}(d)(i + j + k + 1)!}{i!j!k!}  \Big(\dratio\Big)^2\\
    & \quad - (1 + \varepsilon)^2  c_{ijk}(d) \alpha(d) \Big(\dratio\Big)^2\frac{n^2(n - 1)!}{i!j!k!(n - 1 - i - j - k)!} \int_\delta^1 v^{i + j + k + 1}(1 - v)^{n - 1 - i - j - k} \dif  v.
\end{align*}
It is ready to check that the second term in the last display converges to $0$ as $n \to \infty$. We thus have
\begin{align*}
&\limsup_{n \to \infty}\frac{n^2(n - 1)!}{i!j!k!(n - 1 - i - j - k)!} \bE[(V - V_1)^i(V - V_2)^j(V_1 + V_2 - V)^k(1 - V)^{n - 1 - i - j - k}]\\ 
\leq& (1 + \varepsilon)^2\alpha(d)\frac{c_{ijk}(d)(i + j + k + 1)!}{i!j!k!}\Big(\dratio\Big)^2.
\end{align*}
Similarly, we can obtain
\[
\begin{aligned}
&\liminf_{n \to \infty}\frac{n^2(n - 1)!}{i!j!k!(n - 1 - i - j - k)!} \bE[(V - V_1)^i(V - V_2)^j(V_1 + V_2 - V)^k(1 - V)^{n - 1 - i - j - k}]\\ 
\geq& (1 - \varepsilon)^2\alpha(d)\frac{c_{ijk}(d)(i + j + k + 1)!}{i!j!k!}\Big(\dratio\Big)^2.
\end{aligned}
\]
Combining the above two bounds, it then holds true that
\begin{align*}
    & \quad \lim_{n \to \infty} n^2\bE[\nu_1(\cA_M(x))^2 \mid X_1 = x] \\
    & =  \lim_{n \to \infty}\sum_{i + k , j + k \leq M - 1} \frac{n^2(n - 1)!}{i!j!k!(n - 1 - i - j - k)!} \bE\Big[(V - V_1)^i(V - V_2)^j(V_1 + V_2 - V)^k(1 - V)^{n - 1 - i - j - k}\Big] \\
    & = \sum_{i +k, j + k \leq M- 1}\alpha(d)\frac{c_{ijk}(d)(i + j + k + 1)!}{i!j!k!}\Big(\dratio\Big)^2 \\
    & = \alpha(M, d)\Big(\dratio\Big)^2.
\end{align*}

~\\
\textbf{Step 3.} We close the proof by calculating the explicit values of $\alpha(M,d)$ when either $M$ or $d$ is 1. The results for $M=1$ is \citet[Theorem 1]{devroye2017measure}. It remains to prove
\[
\alpha(M,1) = M(2M + 1)/2.
\]
The previous arguments yield
\[
\alpha(M, 1) = \alpha(1)\sum_{i + k \leq M - 1, j + k \leq M - 1}\frac{c_{ijk}(1)(i + j + k + 1)!}{i!j!k!}.
\]
By the form of $c_{ijk}(1)$ in Lemma \ref{lem5.2}, we have
\[
\frac{c_{ijk}(1)(i + j + k + 1)!}{i!j!k!}= \begin{cases}
        1, & \mbox{ if } i = j = k = 0,\\
       \frac{1}{3}, & \mbox{ if } \text{only one of }i, j, k  \text{ is equal to }0,\\
        \frac{2}{3}, & \mbox{ if } \text{only one of }i, j, k  \text{ is not equal to } 0 ,\\
       0, & \mbox{ if } i,j,k \neq 0.
    \end{cases}
\]

It can be shown that, under the restriction $i + k \leq M - 1, j + k \leq M - 1$, there are $2M^2 - 5M + 3$ combinations of $(i, j, k)$ satisfying only one of $i, j, k$ is equal to 0; $3M - 3$ combinations satisfying only one of $i, j, k$ is not equal to 0; and one combination satisfying $(i, j, k) = (0, 0, 0)$. Adding together, we have
\[
\sum_{i + k \leq M - 1, j + k \leq M - 1}\frac{c_{ijk}(1)(i + j + k + 1)!}{i!j!k!} = \frac{1}{3}(2M^2 - 5M + 3) + \frac{2}{3}(3M - 3) + 1 = \frac{1}{3}(2M^2 + M).
\]
Since $\alpha(1) = 3/2$, we obtain
\[
\alpha(M,1) = \frac{3}{2}\cdot \frac{1}{3}(2M^2 + M) = \frac{M(2M + 1)}{2},
\]
and thus completes the proof.

\subsection{Proof of Theorem \ref{thm4.2}}

Think of $\nu_0, \nu_1$ as the conditional distributions of $X \mid W = 0$, $X \mid W = 1$. Let $f_0(\cdot), f_1(\cdot)$ be the corresponding Lebesgue density functions.
By Assumption \ref{assump4.5}, we have both $f_0(\cdot)$ and $f_1(\cdot)$ to be continuous in $\cX$.

~\\
{\bf Step 1.} We first prove the following convergence result for $\bE[V^E]$:
    \begin{align*}
   &\lim_{n\to\infty}\bE[V^E]   \\
   =& \frac{1}{M^2}\bE\Big[\sigma_1^2(X)\Big(\frac{\alpha(M, d)}{e(X)} + (\alpha(M, d) - M^2 - M)e(X) + (2M^2 + M - 2\alpha(M, d))\Big)\Big] \\
    &+ \frac{1}{M^2}\bE\Big[\sigma_0^2(X)\Big(\frac{\alpha(M, d)}{1 - e(X)} + (\alpha(M, d) - M^2 - M)(1 - e(X)) + (2M^2 + M - 2\alpha(M, d))\Big)\Big].
    \end{align*}
Recall that
\[
V^E = \frac{1}{n}\sum_{i = 1}^n \Big(1 + \frac{K_M(i)}{M}\Big)^2\sigma_{W_i}^2(X_i).
\]
Let $p = \bP(W_i = 1)$. Since $(X_i, W_i)_{i = 1}^n$ are independent and identically distributed, we have
\begin{align*}
\bE[V^E] =& \bE\Big[\Big(1 + \frac{K_M(i)}{M}\Big)^2\sigma_{W_i}^2(X_i)\Big] \\
        =& \bE\Big[\Big(1 + \frac{K_M(i)}{M}\Big)^2\sigma_1^2(X_i) \mid W_i = 1\Big]p  + \bE\Big[\Big(1 + \frac{K_M(i)}{M}\Big)^2\sigma_0^2(X_i) \mid W_i = 0\Big](1 - p).
\end{align*}
Fix the treatment indicators $\mathbf{W}=(W_1,\ldots,W_n)$ and covariates in the control group $\{X_j\}_{j:W_j = 0}$. 
It is then immediate that
\[
K_M(i) \mid \mathbf{W}, \{X_j\}_{j:W_j = 0}, W_i = 0 \sim \text{Binomial}\Big(n_1, \nu_1(\cA_M(X_i))\Big).
\]
We accordingly have
\begin{align*}
&\bE\Big[K_M(i) \mid \mathbf{W}, \{X_j\}_{j:W_j = 0}, W_i = 0\Big] = n_1\nu_1(\cA_M(X_i)),\\
&\bE\Big[K_M(i)^2 \mid \mathbf{W}, \{X_j\}_{j:W_j = 0}, W_i = 0\Big] = n_1\nu_1(\cA_M(X_i)) + n_1(n_1 - 1)\nu_1(\cA_M(X_i))^2.
\end{align*}

We first introduce the following lemma.
\begin{lemma} \label{lem5.3} \nb{It holds true that
\[
\bE[n_0 \nu_1(\cA_M(x)) \mid X_i = x,   W_i = 0]\quad  \text{ and } \quad   \bE[n_0^2\nu_1(\cA_M(x))^2 \mid X_i = x,  W_i = 0]
\]
are uniformly bounded for all $n$ and $x \in \cX$.
    }
\end{lemma}

Back to $\bE[V^E]$, we have
\begin{align*}
    &  \bE\Big[\left(1 + \frac{K_M(i)}{M}\right)^2\sigma_0^2(X_i) \mid W_i = 0\Big] \\
   =& \bE\Big[\left(1 + \frac{2}{M}n_1\nu_1(\cA_M(X_i)) + \frac{1}{M^2}(n_1\nu_1(\cA_M(X_i)) + n_1(n_1 - 1)\nu_1(\cA_M(X_i))^2\right)\sigma_0^2(X_i)\mid W_i = 0\Big] \\
   =& \bE\Big[\left(1 + \frac{n_1}{n_0}\Big(\frac{2}{M}+ \frac{1}{M^2}\Big)n_0\nu_1(\cA_M(X_i)) + \frac{n_1(n_1 - 1)}{n_0^2}\frac{1}{M^2}n_0^2\nu_1(\cA_M(X_i))^2\right)\sigma_0^2(X_i)\mid W_i = 0\Big].
\end{align*}

By the law of large numbers, we have 
\[
n_1/n_0 \xrightarrow{a.s.} p/(1 - p). 
\]
\nb{By Lemma \ref{lem5.3} and Assumption \ref{assump4.4},  
\begin{align*}
\bE\Big[n_0 \nu_1(\cA_M(x)) \mid X_i = x, W_i = 0\Big],~  \bE\Big[n_0^2\nu_1(\cA_M(x))^2 \mid X_i = x, W_i = 0\Big],\\
~ {\rm and}~ \sigma_w^2(x) \text{ are all uniformly bounded.} 
\end{align*}
}
Adding together yields
\[
\begin{aligned}
   & \bE\Big[\Big(1 + \frac{n_1}{n_0}(\frac{2}{M}+ \frac{1}{M^2})n_0\nu_1(\cA_M(X_i)) + \frac{n_1(n_1 - 1)}{n_0^2}\frac{1}{M^2}n_0^2\nu_1(\cA_M(X_i))^2\Big)\sigma_0^2(X_i)\mid W_i = 0\Big] \\
   & = \bE\Big[\Big(1 + \frac{p}{1- p}\Big(\frac{2}{M}+ \frac{1}{M^2}\Big)n_0\nu_1(\cA_M(X_i)) + \frac{p^2}{(1 - p)^2}\frac{1}{M^2}n_0^2\nu_1(\cA_M(X_i))^2\Big)\sigma_0^2(X_i)\mid W_i = 0\Big] + o(1).
\end{aligned}
\]
Let
\[
F_n(x) = \bE\Big[\Big(1 + \frac{p}{1- p}\Big(\frac{2}{M}+ \frac{1}{M^2}\Big)n_0\nu_1(\cA_M(X_i)) + \frac{p^2}{(1 - p)^2}\frac{1}{M^2}n_0^2\nu_1(\cA_M(X_i))^2\Big)\mid X_i = x, W_i = 0\Big].
\]
Theorem \ref{thm3.1} then implies, for almost all $x \in \cX$,
\[
F_n(x) \to F(x) = 1 + \frac{p}{1- p}\Big(\frac{2}{M}+ \frac{1}{M^2}\Big)M\cdot\dratio + \frac{p^2}{(1 - p)^2}\frac{1}{M^2}\cdot\alpha(M,d)\Big(\dratio\Big)^2.
\]
\nb{By Lemma \ref{lem5.3}, we have $F_n(x)$ are uniformly bounded for all $n$ and $x \in \cX$. In addition, by Assumption \ref{assump4.3}(ii) and \ref{assump4.4}(iii), we have $F(x)$ are uniformly bounded. Thus 
\[
F_n(x)\sigma_0^2(x)~~~{\rm and}~~~F(x)\sigma_0^2(x) 
\]
are uniformly bounded for all $n$ and $x \in \cX$. Then, by the dominated convergence theorem, we have}
\begin{align*}
&\bE\Big[\Big(1 + \frac{p}{1- p}\Big(\frac{2}{M}+ \frac{1}{M^2}\Big)n_0\nu_1(\cA_M(X_i)) + \frac{p^2}{(1 - p)^2}\frac{1}{M^2}n_0^2\nu_1(\cA_M(X_i))^2\Big)\sigma_0^2(X_i)\mid W_i = 0\Big] \\
 =& \bE\Big[F_n(X_i)\sigma_0^2(X_i) \mid W_i = 0\Big]\\
 =& \bE\Big[F(X_i) \sigma_0^2(X_i) \mid W_i = 0\Big]  + o(1) \\
 =& \bE\Big[\Big(1 + \frac{p}{1- p}\Big(\frac{2}{M}+ \frac{1}{M^2}\Big)M\frac{f_1(X_i)}{f_0(X_i)} + \frac{p^2}{(1 - p)^2}\frac{1}{M^2}\alpha(M,d)\Big(\frac{f_1(X_i)}{f_0(X_i)}\Big)^2\Big)\sigma_0^2(X_i)\mid W_i = 0\Big] + o(1).
\end{align*}
Furthermore, we have
\[
\bE[\sigma_0^2(X_i) \mid W_i = 0](1 - p) =  \bE\Big[\frac{1 - W_i}{1 - p}\sigma_0^2(X_i)\Big](1 - p) =  \bE[\sigma_0^2(X_i)(1 - W_i)] = \bE[\sigma_0^2(X_i)(1 - e(X_i))].
\]
Note  that 
\[
\frac{f_1(X_i)}{f_0(X_i)} = \frac{e(X_i)/p}{(1 - e(X_i))/(1 - p)}.
\]
We accordingly have
\begin{align*}
&\bE\Big[\frac{p}{1 - p}\Big(\frac{2}{M} + \frac{1}{M^2}\Big)M\frac{f_1(X_i)}{f_0(X_i)}\sigma_0^2(X_i) \mid W_i = 0\Big](1 - p)\\
= & \bE\Big[\frac{p}{1 - p}\Big(\frac{2}{M} + \frac{1}{M^2}\Big)M\frac{f_1(X_i)}{f_0(X_i)}\sigma_0^2(X_i) (1 - W_i)\Big]\\
= & \bE\Big[\frac{p}{1 - p}\Big(\frac{2}{M} + \frac{1}{M^2}\Big)M\frac{f_1(X_i)}{f_0(X_i)}\sigma_0^2(X_i) (1 - e(X_i))\Big]\\
= & \bE\Big[\Big(\frac{2}{M} + \frac{1}{M^2}\Big)M\frac{e(X_i)}{1 - e(X_i)}\sigma_0^2(X_i) (1 - e(X_i))\Big] \\
= & \bE\Big[(2 + \frac{1}{M})e(X_i)\sigma_0^2(X_i)\Big],
\end{align*}
and similarly,
\begin{align*}
&\bE\Big[\frac{p^2}{(1 - p)^2}\frac{1}{M^2}\alpha(M,d)\Big(\frac{f_1(X_i)}{f_0(X_i)}\Big)^2 \sigma_0^2(X_i)\mid W_i = 0\Big](1 - p) \\
= & \bE\Big[\frac{p^2}{(1 - p)^2}\frac{1}{M^2}\alpha(M,d)\Big(\frac{f_1(X_i)}{f_0(X_i)}\Big)^2 \sigma_0^2(X_i)(1 - W_i)\Big] \\
= & \bE\Big[\frac{p^2}{(1 - p)^2}\frac{1}{M^2}\alpha(M,d)\Big(\frac{f_1(X_i)}{f_0(X_i)}\Big)^2 \sigma_0^2(X_i)(1 - e(X_i))\Big] \\
= &  \bE\Big[\frac{\alpha(M,d)}{M^2}\frac{e(X_i)^2}{(1 - e(X_i))^2}\sigma_0^2(X_i)(1 - e(X_i))\Big]\\
= & \bE\Big[\frac{\alpha(M,d)}{M^2}\frac{e(X_i)^2}{1 - e(X_i)}\sigma_0^2(X_i)\Big].
\end{align*}

Combining the above two identities yields
\begin{align*}
  & \bE\Big[\Big(1 + \frac{K_M(i)}{M}\Big)^2\sigma_0^2(X_i) \mid W_i = 0\Big](1 - p)\\
   = &\bE\Big[\Big(1 + \frac{p}{1- p}\Big(\frac{2}{M}+ \frac{1}{M^2}\Big)M\frac{f_1(X_i)}{f_0(X_i)} + \frac{p^2}{(1 - p)^2}\frac{1}{M^2}\alpha(M,d)\Big(\frac{f_1(X_i)}{f_0(X_i)}\Big)^2\Big)\sigma_0^2(X_i)\mid W_i = 0\Big] + o(1) \\
    = &  \frac{1}{M^2}\bE\Big[\sigma_0^2(X)\Big(\frac{\alpha(M, d)}{1 - e(X)} + (\alpha(M, d) - M^2 - M)(1 - e(X)) + (2M^2 + M - 2\alpha(M, d))\Big)\Big] + o(1).
\end{align*}

Similarly, we have
\begin{align*}
  & \bE\Big[\Big(1 + \frac{K_M(i)}{M}\Big)^2\sigma_1^2(X_i) \mid W_i = 1\Big]p\\
    = &  \frac{1}{M^2}\bE\Big[\sigma_1^2(X)\Big(\frac{\alpha(M, d)}{e(X)} + (\alpha(M, d) - M^2 - M)e(X) + (2M^2 + M - 2\alpha(M, d))\Big)\Big] + o(1).
\end{align*}

Combine these two yields the convergence result for $\bE[V^E]$.

\vspace{0.5cm}

\nb{
{\bf Step 2.} In this step, we show that $\var(V^E)\to 0$, which implies $V^E - \bE[V^E]\xrightarrow{\enskip\P\enskip} 0$.

Following the argument of \citet[Page 40]{abadie2002simple}, it suffices to establish that, for any $i\neq j$,
\begin{align*}\label{eq:negcov}
    \limsup_{n\to\infty}\Big(&\bE\Big[\Big(1 + \frac{K_M(i)}{M}\Big)^2\Big(1 + \frac{K_M(j)}{M}\Big)^2\sigma^2_{W_i}(X_i)\sigma^2_{W_j}(X_j)\Big]\\
    &-\bE\Big[\Big(1 + \frac{K_M(i)}{M}\Big)^2\sigma^2_{W_i}(X_i)\Big]\bE\Big[\Big(1 + \frac{K_M(j)}{M}\Big)^2\sigma^2_{W_j}(X_j)\Big]\Big)\leq 0.
\end{align*}

{\bf Case 1: Same arm.} Consider the case where $W_i=W_j=0$; the case $W_i=W_j=1$ follows by symmetry.

For each control unit $k$ with $W_k=0$, define
\[
p_k:=\nu_1(\cA_M(X_k))=\int_{\cA_M(X_k)}f_1(z)\dif z, \qquad q_k:=\int_{\cA_M(X_k)}f_1(z)\sigma_0^2(z)\dif z.
\]
For $i\neq j$ with $W_i=W_j=0$, further define
\[
p_{ij}:=\nu_1\Big(\cA_M(X_i)\cap\cA_M(X_j)\Big).
\]
Conditioning on $\mathbf{W}$ and $\{X_k\}_{W_k=0}$, and using that $\{X_k\}_{W_k=1}$ are i.i.d.\ given $\mathbf{W}$, we have
\begin{align*}
&\bE\Big[K_M(i)\sigma_0^2(X_i)\;\Big|\;\mathbf{W},\{X_k\}_{k:W_k=0},W_i=0\Big]=n_1\,p_i\,\sigma_0^2(X_i),\\
&\bE\Big[K_M(i)K_M(j)\sigma_0^2(X_i)\sigma_0^2(X_j)\;\Big|\;\mathbf{W},\{X_k\}_{k:W_k=0},W_i=W_j=0\Big]\\
&\qquad =\big(n_1p_{ij}+n_1(n_1-1)p_ip_j\big)\sigma_0^2(X_i)\sigma_0^2(X_j).
\end{align*}

We now relate $p_k\sigma_0^2(X_k)$ to $q_k$ via the Lipschitz property of $\sigma_0^2(\cdot)$. By Assumption \ref{assump4.4}(i), there exists a constant $L>0$ such that, for any $z\in\cA_M(X_k)$,
\[
\big|\sigma_0^2(X_k)-\sigma_0^2(z)\big|\leq L\,\|X_k-z\|\leq L\cdot\mathrm{diam}\big(\cA_M(X_k)\big),
\]
where $\mathrm{diam}$ represents the diameter of the input set measured in $\|\cdot\|$.

Consequently,
\begin{align}\label{eq:lip-pq}
\big|p_k\sigma_0^2(X_k)-q_k\big|=\Big|\int_{\cA_M(X_k)}f_1(z)\big(\sigma_0^2(X_k)-\sigma_0^2(z)\big)\dif z\Big|\leq L\,p_k\cdot\mathrm{diam}\big(\cA_M(X_k)\big).
\end{align}

To bound this error in \eqref{eq:lip-pq}, we introduce the following lemma.

\begin{lemma}\label{lem5.4}
Under Assumptions \ref{assump4.1}--\ref{assump4.5}, for any fixed $i$,
\[
\mathrm{diam}\big(\cA_M(X_i)\big)\xrightarrow{\enskip\bP\enskip}0,\qquad\text{as }n\to\infty.
\]
Further more, $\bE[\mathrm{diam}\big(\cA_M(X_i)\big)] \to 0$ as $n \to \infty$.
\end{lemma}

By definition of the catchment area, each treated unit belongs to the catchment areas of exactly $M$ control units, so
\begin{align}\label{eq:han-1}
\sum_{k:\,W_k=0}q_k=\int f_1(z)\sigma_0^2(z)\sum_{k:\,W_k=0}\ind\big\{z\in\cA_M(X_k)\big\}\dif z=M\int f_1(z)\sigma_0^2(z)\dif z,
\end{align}
which is a deterministic constant. By exchangeability of $\{X_k\}_{W_k=0}$ given $\mathbf{W}$ and the exact constraint \eqref{eq:han-1},
\[
\bE\big[q_i\;\big|\;\mathbf{W},W_i=0\big]=\frac{1}{n_0}\sum_{k:\,W_k=0}\bE\big[q_k\;\big|\;\mathbf{W}\big]=\frac{M}{n_0}\int f_1(z)\sigma_0^2(z)\dif z.
\]
Combined with \eqref{eq:lip-pq}, the boundedness of $n_0p_i$ in $L^2$ (Lemma \ref{lem5.3}), Lemma \ref{lem5.4}, and the dominated convergence theorem, we obtain
\[
\bE\big[n_1\big|p_i\sigma_0^2(X_i)-q_i\big|\;\big|\;W_i=0\big]\leq L\,\bE\big[n_1p_i\cdot\mathrm{diam}(\cA_M(X_i))\;\big|\;W_i=0\big]\to 0.
\]
Combining the above, we obtain the marginal limit
\begin{align}\label{eq:marginal-K}
\bE\big[K_M(i)\sigma_0^2(X_i)\;\big|\;W_i=0\big]=\bE\big[n_1q_i\;\big|\;W_i=0\big]+o(1)\to\frac{p}{1-p}\,M\int f_1(z)\sigma_0^2(z)\dif z.
\end{align}

For the cross term, by exchangeability and the constraint \eqref{eq:han-1}, we have
\begin{align*}
\bE\Big[q_iq_j\;\Big|\;\mathbf{W},W_i=W_j=0\Big]&=\bE\Big[q_i\,\frac{\sum_{k\neq i,W_k=0}q_k}{n_0-1}\;\Big|\;\mathbf{W},W_i=W_j=0\Big]\\
&=\bE\Big[q_i\,\frac{M\int f_1\sigma_0^2-q_i}{n_0-1}\;\Big|\;\mathbf{W},W_i=W_j=0\Big]\\
&\leq\frac{M\int f_1\sigma_0^2}{n_0-1}\,\bE\big[q_i\;\big|\;\mathbf{W},W_i=W_j=0\big]\\
&=\frac{M^2\big(\int f_1\sigma_0^2\big)^2}{n_0(n_0-1)}.
\end{align*}
Replacing $p_ip_j\sigma_0^2(X_i)\sigma_0^2(X_j)$ by $q_iq_j$ creates an error bounded, via \eqref{eq:lip-pq} and Assumption \ref{assump4.4}(iii), by
\begin{equation}\label{error}
\big|p_ip_j\sigma_0^2(X_i)\sigma_0^2(X_j)-q_iq_j\big|\leq L\overline{\sigma}^2\,p_ip_j\,\big(\mathrm{diam}(\cA_M(X_i))+\mathrm{diam}(\cA_M(X_j))\big).
\end{equation}
The righthand side of \eqref{error} vanishes upon multiplying by $n_1(n_1-1)$ and taking expectation by Cauchy--Schwarz, Lemmas \ref{lem5.3} and \ref{lem5.4}, and dominated convergence. 

For the remaining $n_1p_{ij}$ term, we have
\[
\bE\big[n_1p_{ij}\sigma_0^2(X_i)\sigma_0^2(X_j)\;\big|\;W_i=W_j=0\big]\leq\overline{\sigma}^4\,\bE\big[n_1p_{ij}\;\big|\;W_i=W_j=0\big]\to 0,
\]
where the last convergence follows since 
\[
\bE[n_1p_{ij}\mid\mathbf{W},W_i=W_j=0]\leq Mn_1/(n_0(n_0-1))\to 0 
\]
by the identity $\sum_{k:W_k=0}\ind\{z\in\cA_M(X_k)\}\leq M$.

Combining the above, we obtain
\begin{align}\label{eq:cross-K-same}
\limsup_{n \to \infty}\bE\Big[K_M(i)K_M(j)\sigma_0^2(X_i)\sigma_0^2(X_j)\;\Big|\;W_i=W_j=0\Big]\leq\frac{p^2}{(1-p)^2}M^2\Big(\int f_1(z)\sigma_0^2(z)\dif z\Big)^2,
\end{align}
and the righthand side exactly matches the product of marginal limits in \eqref{eq:marginal-K}. Hence
\begin{align*}
\limsup_{n \to \infty}\Big(&\bE\Big[K_M(i)K_M(j)\sigma_0^2(X_i)\sigma_0^2(X_j)\;\Big|\;W_i=W_j=0\Big]\\
&-\bE\big[K_M(i)\sigma_0^2(X_i)\,\big|\,W_i=0\big]\bE\big[K_M(j)\sigma_0^2(X_j)\,\big|\,W_j=0\big]\Big)\leq   0.
\end{align*}
The case $W_i=W_j=1$ is analogous, with $\sigma_0^2,f_1,p/(1-p)$ replaced by $\sigma_1^2,f_0,(1-p)/p$, respectively.

\vspace{0.3cm}

{\bf Case 2: Opposite arm.} Consider the case where $W_i=0$ and $W_j=1$; the case $W_i=1$, $W_j=0$ follows by symmetry. 

We adapt the argument used in the unweighted case by attaching $\sigma_0^2(X_i)$ to each indicator $\ind\{X_k\in\cA_M(X_i)\}$ and $\sigma_1^2(X_j)$ to each $\ind\{X_\ell\in\cA_M(X_j)\}$. Recall that
\begin{align*}
K_M(i)\sigma_0^2(X_i)=\sum_{k:\,W_k=1}\ind\big\{X_k\in\cA_M(X_i)\big\}\sigma_0^2(X_i),\\ 
K_M(j)\sigma_1^2(X_j)=\sum_{\ell:\,W_\ell=0}\ind\big\{X_\ell\in\cA_M(X_j)\big\}\sigma_1^2(X_j),
\end{align*}
where we recall that 
\[
\cA_M(X_i)=\cA_M(X_i;\{X_k\}_{W_k=0}) 
\]
is defined with respect to the control sample and 
\[
\cA_M(X_j)=\cA_M(X_j;\{X_k\}_{W_k=1}) 
\]
with respect to the treated sample. For each pair $(k,\ell)$ with $W_k=1$ and $W_\ell=0$, we then define the reduced catchment areas
\begin{align*}
\cA_{M,-\ell}(X_i):=\cA_M\big(X_i;\{X_m\}_{W_m=0}\setminus\{X_\ell\}\big),\\
\cA_{M,-k}(X_j):=\cA_M\big(X_j;\{X_m\}_{W_m=1}\setminus\{X_k\}\big).
\end{align*}
Since removing a unit from the reference sample can only enlarge the catchment area, we have the following the key observation that 
\[
\cA_M(X_i)\subseteq\cA_{M,-\ell}(X_i) ~~{\rm and}~~ \cA_M(X_j)\subseteq\cA_{M,-k}(X_j), 
\]
so that
\begin{align*}
\ind\big\{X_k\in\cA_M(X_i)\big\}\leq\ind\big\{X_k\in\cA_{M,-\ell}(X_i)\big\},\\
\ind\big\{X_\ell\in\cA_M(X_j)\big\}\leq\ind\big\{X_\ell\in\cA_{M,-k}(X_j)\big\}.
\end{align*}
Pick a specific pair $(i_0,j_0)$ such that $i_0\neq i$, $j_0\neq j$ and $W_{i_0}=0$, $W_{j_0}=1$. By exchangeability of $\{X_k\}_{W_k=0}$ and of $\{X_\ell\}_{W_\ell=1}$ given $\mathbf{W}$, we have
\begin{align*}
&\quad \bE\Big[K_M(i)K_M(j)\sigma_0^2(X_i)\sigma_1^2(X_j)\;\Big|\;\mathbf{W},W_i=0,W_j=1\Big] \\
& =\bE\Big[\Big(\sum_{k:\,W_k=1}\ind\{X_k\in\cA_M(X_i)\}\Big)\Big(\sum_{\ell:\,W_\ell=0}\ind\{X_\ell\in\cA_M(X_j)\}\Big)
\sigma_0^2(X_i)\sigma_1^2(X_j)\;\Big|\;\mathbf{W},W_i=0,W_j=1\Big] \\
& =\bE\Big[n_0n_1\,\ind\{X_{j_0}\in\cA_M(X_i)\}\ind\{X_{i_0}\in\cA_M(X_j)\}\sigma_0^2(X_i)\sigma_1^2(X_j)\;\Big|\;\mathbf{W},W_i=0,W_j=1\Big] \\
& \leq\bE\Big[n_0n_1\,\ind\{X_{j_0}\in\cA_{M,-i_0}(X_i)\}\ind\{X_{i_0}\in\cA_{M,-j_0}(X_j)\}\sigma_0^2(X_i)\sigma_1^2(X_j)\;\Big|\;\mathbf{W},W_i=0,W_j=1\Big].
\end{align*}
Conditioning on $\mathbf{W}$, the random object 
\[
\ind\{X_{j_0}\in\cA_{M,-i_0}(X_i)\}\sigma_0^2(X_i) 
\]
depends only through $X_{j_0}$ and $\{X_m\}_{W_m=0}\setminus\{X_{i_0}\}$, while 
\[
\ind\{X_{i_0}\in\cA_{M,-j_0}(X_j)\}\sigma_1^2(X_j) 
\]
depends only through $X_{i_0}$ and $\{X_\ell\}_{W_\ell=1}\setminus\{X_{j_0}\}$, and these two collections are independent given $\mathbf{W}$. Therefore the righthand side factorizes, yielding
\begin{align*}
&\bE\Big[K_M(i)K_M(j)\sigma_0^2(X_i)\sigma_1^2(X_j)\;\Big|\;\mathbf{W},W_i=0,W_j=1\Big]\\ 
&\qquad\leq\bE\big[n_1\,\ind\{X_{j_0}\in\cA_{M,-i_0}(X_i)\}\sigma_0^2(X_i)\;\big|\;\mathbf{W},W_i=W_{i_0}=0\big]\\
&\qquad\qquad\times\bE\big[n_0\,\ind\{X_{i_0}\in\cA_{M,-j_0}(X_j)\}\sigma_1^2(X_j)\;\big|\;\mathbf{W},W_j=W_{j_0}=1\big].
\end{align*}

To handle each factor, define the reduced-area analogue of $q_k$:
\[
q_k^{(-i_0)}:=\int_{\cA_{M,-i_0}(X_k)}f_1(z)\sigma_0^2(z)\dif z,\qquad k\in\Big\{m:\,W_m=0,\,m\neq i_0\Big\}.
\]
Marginalizing over $X_{j_0}$ first and applying the Lipschitz argument \eqref{eq:lip-pq} to the reduced catchment area, we have
\begin{align*}
&\bE\big[\ind\{X_{j_0}\in\cA_{M,-i_0}(X_i)\}\sigma_0^2(X_i)\;\big|\;\mathbf{W},\{X_m\}_{W_m=0}\setminus\{X_{i_0}\},W_i=W_{i_0}=0\big]\\
=&\nu_1(\cA_{M,-i_0}(X_i))\,\sigma_0^2(X_i)=:q_i^{(-i_0)}+\eta_i,
\end{align*}
where the second term above satisfies
\[
|\eta_i|\leq L\,\nu_1(\cA_{M,-i_0}(X_i))\cdot\mathrm{diam}(\cA_{M,-i_0}(X_i)). 
\]

Since each treated unit belongs to exactly $M$ of the reduced catchment areas 
\[
\Big\{\cA_{M,-i_0}(X_k)\Big\}_{k\neq i_0,W_k=0}, 
\]
we have
 \[
 \sum_{k\neq i_0,W_k=0}q_k^{(-i_0)}=M\int f_1\sigma_0^2, 
 \]
 and exchangeability of $\{X_k\}_{k\neq i_0,W_k=0}$ given $\mathbf{W}$ gives
\[
\bE\big[q_i^{(-i_0)}\;\big|\;\mathbf{W},W_i=W_{i_0}=0\big]=\frac{M}{n_0-1}\int f_1(z)\sigma_0^2(z)\dif z.
\]
Combined with Lemma \ref{lem5.4} and dominated convergence, marginalizing over $\mathbf{W}$ yields
\[
\bE\big[n_1\,\ind\{X_{j_0}\in\cA_{M,-i_0}(X_i)\}\sigma_0^2(X_i)\;\big|\;W_i=0,W_j=1\big]\to\frac{p}{1-p}\,M\int f_1(z)\sigma_0^2(z)\dif z.
\]
By a symmetric argument with $\sigma_0^2,f_1$ replaced by $\sigma_1^2,f_0$,
\[
\bE\big[n_0\,\ind\{X_{i_0}\in\cA_{M,-j_0}(X_j)\}\sigma_1^2(X_j)\;\big|\;W_i=0,W_j=1\big]\to\frac{1-p}{p}\,M\int f_0(z)\sigma_1^2(z)\dif z.
\]
Combining, we obtain
\[
\limsup_{n\to\infty}\bE\Big[K_M(i)K_M(j)\sigma_0^2(X_i)\sigma_1^2(X_j)\;\Big|\;W_i=0,W_j=1\Big]\leq M^2\int f_1(z)\sigma_0^2(z)\dif z\int f_0(z)\sigma_1^2(z)\dif z.
\]
On the other hand, by the marginal limit \eqref{eq:marginal-K} and its analogue for $W_j=1$,
\begin{align*}
\bE\big[K_M(i)\sigma_0^2(X_i)\,\big|\,W_i=0\big]\bE\big[K_M(j)\sigma_1^2(X_j)\,\big|\,W_j=1\big]\to\frac{p}{1-p}M\int f_1\sigma_0^2\cdot\frac{1-p}{p}M\int f_0\sigma_1^2\\
=M^2\int f_1\sigma_0^2\int f_0\sigma_1^2.
\end{align*}
Hence
\begin{align*}
\limsup_{n \to \infty}\Big(&\bE\Big[K_M(i)K_M(j)\sigma_0^2(X_i)\sigma_0^2(X_j)\;\Big|\;W_i=W_j=0\Big]\\
&-\bE\big[K_M(i)\sigma_0^2(X_i)\,\big|\,W_i=0\big]\bE\big[K_M(j)\sigma_0^2(X_j)\,\big|\,W_j=0\big]\Big)\leq   0.
\end{align*}

\vspace{0.2cm}

Combining {\bf Cases 1 and 2}, we obtain, for any $i\neq j$,
\[
\limsup_{n\to\infty}\Big(\bE\big[K_M(i)K_M(j)\sigma^2_{W_i}(X_i)\sigma^2_{W_j}(X_j)\big]-\bE\big[K_M(i)\sigma^2_{W_i}(X_i)\big]\bE\big[K_M(j)\sigma^2_{W_j}(X_j)\big]\Big)\leq 0.
\]

Note that exactly the same proof technique also yields
\[
\limsup_{n\to\infty}\Big(\bE\big[K_M(i)^2K_M(j)\sigma^2_{W_i}(X_i)\sigma^2_{W_j}(X_j)\big]-\bE\big[K_M(i)^2\sigma^2_{W_i}(X_i)\big]\bE\big[K_M(j)\sigma^2_{W_j}(X_j)\big]\Big)\leq 0
\]
and
\[
\limsup_{n\to\infty}\Big(\bE\big[K_M(i)^2K_M(j)^2\sigma^2_{W_i}(X_i)\sigma^2_{W_j}(X_j)\big]-\bE\big[K_M(i)^2\sigma^2_{W_i}(X_i)\big]\bE\big[K_M(j)^2\sigma^2_{W_j}(X_j)\big]\Big)\leq 0.
\]
By multiplying different constant and adding up for the above results, we then obtain
\begin{align*}\label{eq:negcov}
    \limsup_{n\to\infty}\Big(&\bE\Big[\Big(1 + \frac{K_M(i)}{M}\Big)^2\Big(1 + \frac{K_M(i)}{M}\Big)^2\sigma^2_{W_i}(X_i)\sigma^2_{W_j}(X_j)\Big]\\
    &-\bE\Big[\Big(1 + \frac{K_M(i)}{M}\Big)^2\sigma^2_{W_i}(X_i)\Big]\bE\Big[\Big(1 + \frac{K_M(i)}{M}\Big)^2\sigma^2_{W_j}(X_j)\Big]\Big)\leq 0,
\end{align*}
which proves $\Var(V^E) \to 0$ as $n \to \infty$.

~\\
{\bf Step 3.} Finally, we establish asymptotic normality using the martingale representation of the matching estimator, following the argument of \citet{abadie2012martingale}, which was also exploited by \citet{lin2023estimation}. This provides an alternative proof to that of \citet{abadie2006large}, whose argument contains a logical gap that was later addressed in \citet{abadie2012martingale}. We include the proof here mainly for technical completeness.

In detail, define a two-stage filtration $\{\cF_{n,k}\}_{k=0}^{2n}$ by
\begin{align*}
\cF_{n,k} := \sigma\{W_1,\ldots,W_n,\,X_1,\ldots,X_k\},\quad 0\leq k\leq n,\\
\cF_{n,n+k} := \cF_{n,n}\vee\sigma\{Y_1,\ldots,Y_k\},\quad 1\leq k\leq n,
\end{align*}
and define the array $\{\xi_{n,k},\cF_{n,k}\}_{k=1}^{2n}$ of differences by
\[
\xi_{n,k}:=\begin{cases}
\dfrac{1}{\sqrt{n}}\bigl(\tau(X_k)-\tau\bigr), & 1\leq k\leq n,\\[6pt]
\dfrac{1}{\sqrt{n}}(2W_{k-n}-1)\Bigl(1+\dfrac{K_M(k-n)}{M}\Bigr)\varepsilon_{k-n}, & n+1\leq k\leq 2n.
\end{cases}
\]
Since $X_k$ is independent of $\cF_{n,k-1}$ with $\bE[\tau(X_k)]=\tau$, and that $\bE[\varepsilon_i\mid \mathbf{W},\mathbf{X}]=0$, we have 
\[
\bE[\xi_{n,k}\mid\cF_{n,k-1}]=0, 
\]
so that $\{\xi_{n,k},\cF_{n,k}\}$ is a martingale difference array. Moreover, using the decomposition 
\[
\hat{\tau}_M - \tau = \overline{\tau(X)} + E_M + B_M - \tau, 
\]
a direct calculation gives
\[
\sum_{k=1}^{2n}\xi_{n,k} = \sqrt{n}(\hat{\tau}_M - B_M - \tau).
\]
This shows that 
\[
\sum_{k=1}^{2n}\bE[\xi_{n,k}^2\mid\cF_{n,k-1}]\xrightarrow{\;\bP\;}\sigma_{M,d}^2.
\]

For the first stage, since $X_k$ is independent of $\cF_{n,k-1}$, we have
\[
\sum_{k=1}^{n}\bE\bigl[\xi_{n,k}^2\mid\cF_{n,k-1}\bigr]
=\frac{1}{n}\sum_{k=1}^{n}\var\bigl(\tau(X)\bigr) = V^{\tau(X)}.
\]
For the second stage, since $\bE[\varepsilon_i^2\mid W_i=w,X_i=x]=\sigma_w^2(x)$ and $Y_i$ is independent of $Y_j$ given $(\mathbf{W},\mathbf{X})$ for $i\neq j$, we have
\[
\sum_{k=n+1}^{2n}\bE\bigl[\xi_{n,k}^2\mid\cF_{n,k-1}\bigr]
=\frac{1}{n}\sum_{i=1}^{n}\Bigl(1+\frac{K_M(i)}{M}\Bigr)^2\sigma_{W_i}^2(X_i) = V^E.
\]
By {\bf Step~2}, 
\[
V^E-\bE[V^E]\xrightarrow{\;\bP\;}0;
\]
by {\bf Step~1}, 
\[
\bE[V^E]\to\sigma_{M,d}^2-V^{\tau(X)}. 
\]
Therefore
\[
\sum_{k=1}^{2n}\bE\bigl[\xi_{n,k}^2\mid\cF_{n,k-1}\bigr] = V_M \xrightarrow{\;\bP\;} \sigma_{M,d}^2.
\]
To apply a Martingale Central Limit Theorem, it suffices to check the Lindeberg condition, that is, for any $\epsilon > 0$, 
\[
\sum_{k=1}^{2n}\bE\bigl[\xi_{n,k}^2\ind\{|\xi_{n,k}| > \epsilon\}\bigr] \to 0,
\]
which is implied by the following Lyapunov condition:
\[
\bE\Bigl[\sum_{k=1}^{2n}|\xi_{n,k}|^{2 + \delta}\Bigr] \to 0, 
\]
for some $\delta > 0$.
We verify the Lyapunov condition with $\delta=2$. For the first stage, since $\tau(\cdot)$ is bounded on the compact $\cX$ by Assumption~\ref{assump4.4}(i):
\[
\bE\Bigl[\sum_{k=1}^{n}|\xi_{n,k}|^4\Bigr]
=\frac{1}{n^2}\sum_{k=1}^{n}\bE\bigl[(\tau(X_k)-\tau)^4\bigr] = O(n^{-1})\to 0.
\]
For the second stage, Assumption~\ref{assump4.4}(ii) gives $\bE[\varepsilon_i^4\mid W_i,X_i]\leq C$ uniformly. A standard moment calculation for the conditional binomial distribution of $K_M(i)$, analogous to Lemma~\ref{lem5.3}, shows that $\bE[(1+K_M(i)/M)^4]$ is uniformly bounded in $n$. Therefore
\[
\bE\Bigl[\sum_{k=n+1}^{2n}|\xi_{n,k}|^4\Bigr]
\leq \frac{C}{n^2}\sum_{i=1}^{n}\bE\Bigl[\Bigl(1+\frac{K_M(i)}{M}\Bigr)^4\Bigr] = O(n^{-1})\to 0.
\]
By Markov's inequality, 
\[
\sum_{k=1}^{2n}\bE[|\xi_{n,k}|^4\mid\cF_{n,k-1}]\xrightarrow{\;\bP\;}0, 
\]
which implies the Lindeberg condition.

Applying the martingale central limit theorem \citep{MR1324786}, we conclude
\[
\sqrt{n}(\hat{\tau}_M - B_M - \tau) \text{ converges in distribution to }\cN(0,\sigma_{M,d}^2),
\]
proving part~(i). Part~(ii) follows from $\sqrt{n}(\hat{B}_M-B_M)\xrightarrow{\;\bP\;}0$ and Slutsky's theorem, which completes the whole proof.
}

\section{Proofs of the rest results}

\subsection{Proof of Lemma \ref{lem5.1}}

The first and second parts correspond to the proofs of Equation (2.2) and the fourth identity in Part (ii) of Theorem 2.1 in \cite{devroye2017measure}, respectively. The proofs then only take minor modifications with regard to the change of measures, and are accordingly omitted.

\subsection{Proof of Lemma \ref{lem5.2}}
\textbf{Step 1}: We first consider the case of dimension $d \geq 2$.

Denote the conditional density functions of $V \mid V_1, V_2$ and $V_1, V_2 \mid V$ by $f_{V \mid V_1, V_2}(\cdot \mid \cdot, \cdot)$ and $f_{V_1, V_2 \mid V_2}(\cdot ,\cdot \mid  \cdot)$, respectively. Let 
\[
D(v) = \Big\{(v_1, v_2) \in \bR^2: 0 \leq v_1, v_2 \leq v, v_1 + v_2 \geq v\Big\}
\]
be the support of $V_1, V_2 \mid V = v$. Thus we have
\begin{align*}
    Q_{ijk}(v) &= \bE\Big[P_{ijk}(v, V_1, V_2) \mid V = v\Big] \\
    & = \int_{D(v)}P_{ijk}(v,v_1,v_2) f_{V_1, V_2\mid v}(v_1, v_2\mid v) \dif v_1 \dif v_2 \\
    & = \int_{D(v)}P_{ijk}(v,v_1,v_2) \frac{f_{V \mid V_1, V_2}(v\mid v_1, v_2) f_{V_1, V_2}(v_1, v_2)}{f_V(v)} \dif v_1 \dif v_2 \\
    & = \int_{D(1)}P_{ijk}(v, vv_1, vv_2)\frac{f_{V \mid V_1, V_2}(v\mid vv_1, vv_2) f_{V_1, V_2}(vv_1, vv_2)}{f_V(v)}v^2 \dif v_1 \dif v_2\\
    & = \int_{D(1)}P_{ijk}(1, v_1, v_2)(vf_{V \mid V_1, V_2}(v\mid vv_1, vv_2))\frac{f_{V_1, V_2}(vv_1, vv_2)}{f_V(v)/v} \dif v_1 \dif v_2.
\end{align*}
By Lemma \ref{lem5.1}, for fixed $v_1, v_2$, as $v \to 0+$, we have
\[
f_{V_1, V_2}(vv_1, vv_2) = f_{V_1}(vv_1)f_{V_2}(vv_2) \to \Big(\dratio\Big)^2
~~~{\rm and}~~~
\frac{f_V(v)}{v} \to \alpha(d)\Big(\dratio\Big)^2.
\]
It suffice to consider the convergence of $vf_{V \mid V_1, V_2}(v\mid vv_1, vv_2)$.

Introduce
\[
R_1 = \Vert Z_1 - x\Vert,~~ \theta_1 =\frac{Z_1 - x}{\Vert Z_1 - x\Vert},~~ R_2 = \Vert Z_2 - x\Vert,~~ \theta_2 = \frac{Z_2 - x}{\Vert Z_2 - x \Vert}. 
\]
Note that the volume of $B(Z_1, \Vert Z_1 - x\Vert) \cup B(Z_2, \Vert Z_2 - x\Vert) $  is uniquely determined by $R_1, R_2$, and the directions $\theta_1, \theta_2$. We thus denote this volume by
\[
\lambda\Big(B(Z_1, \Vert Z_1 - x\Vert) \cup B(Z_2, \Vert Z_2 - x\Vert)\Big) =: S(R_1, R_2, \theta_1, \theta_2).
\]
Then for any $k > 0$, we have
\[
S(kR_1, kR_2, \theta_1, \theta_2) = k^dS(R_1, R_2, \theta_1, \theta_2).
\]
By the argument in Lemma \ref{lem5.1}, for any $\varepsilon > 0$, there exists $\delta > 0$ such that if $V < \delta $ holds, then
\[
(1 - \varepsilon)f_0(x) \leq \frac{\nu_0(B(Z_1, \Vert Z_1 - x\Vert) \cup B(Z_2, \Vert Z_2 - x\Vert))}{\lambda(B(Z_1, \Vert Z_1 - x\Vert) \cup B(Z_2, \Vert Z_2 - x\Vert))} = \frac{V}{S(R_1, R_2, \theta_1, \theta_2)} \leq (1 + \varepsilon)f_0(x).
\]
Similarly, letting $A(r) = c_d r^d$ denote the Lebesgue measure of a ball with radius $r$ in $\bR^d$, we have
\[
(1 - \varepsilon)f_0(x) \leq \frac{\nu_0(B(Z_1, \Vert Z_1 - x\Vert) )}{\lambda(B(Z_1, \Vert Z_1 - x\Vert) )} = \frac{V_1}{A(R_1)} \leq (1 + \varepsilon)f_0(x),
\]
\[
(1 - \varepsilon)f_0(x) \leq \frac{\nu_0(B(Z_2, \Vert Z_2 - x\Vert) )}{\lambda(B(Z_2, \Vert Z_2 - x\Vert) )} = \frac{V_2}{A(R_2)} \leq (1 + \varepsilon)f_0(x).
\]
For any fixed $v_1, v_2, t > 0$ and $ 0 < v < \delta$, consider the conditional probability $\bP(V \leq vt \mid V_1 = vv_1, V_2 = vv_2)$. We have
\begin{align*}
&\quad \bP(V \leq vt \mid V_1 = vv_1, V_2 = vv_2) \\
& \leq \bP\Big( (1 - \varepsilon)f_0(x)S(R_1, R_2, \theta_1, \theta_2) \leq vt \mid V_1 = vv_1, V_2 = vv_2\Big) \\
& \leq \bP\Big( (1 - \varepsilon)f_0(x) S\Big(\Big(\frac{V_1}{c_d(1 + \varepsilon)f_0(x)}\Big)^{1/d}, \Big(\frac{V_2}{c_d(1 + \varepsilon)f_0(x)}\Big)^{1/d}, \theta_1, \theta_2  \Big) \leq vt \mid V_1 = vv_1, V_2 = vv_2\Big) \\
& = \bP\Big(\frac{1 - \varepsilon}{1 + \varepsilon}v c_d^{-1}S(v_1^{1/d}, v_2^{1/d}, \theta_1, \theta_2) \leq vt \mid V_1 = vv_1, V_2 = vv_2\Big) \\
& = \bP\Big(c_d^{-1}S(v_1^{1/d}, v_2^{1/d}, \theta_1, \theta_2) \leq \frac{1 + \varepsilon}{1 - \varepsilon}t \mid V_1 = vv_1, V_2 = vv_2\Big). 
\end{align*}
This shows
\begin{align}\label{eq:6.1}
\bP(V \leq vt \mid V_1 = vv_1, V_2 = vv_2) \leq \bP\Big(c_d^{-1}S(v_1^{1/d}, v_2^{1/d}, \theta_1, \theta_2) \leq \frac{1 + \varepsilon}{1 - \varepsilon}t \mid V_1 = vv_1, V_2 = vv_2\Big). 
\end{align}
Similarly, we have 
\begin{align}\label{eq:6.2}
\bP(V \leq vt \mid V_1 = vv_1, V_2 = vv_2) \geq \bP\Big(c_d^{-1}S(v_1^{1/d}, v_2^{1/d}, \theta_1, \theta_2) \leq \frac{1 - \varepsilon}{1 + \varepsilon}t \mid V_1 = vv_1, V_2 = vv_2\Big).
\end{align}
We then consider the limiting distribution of $\theta_1 \mid V_1 = vv_1$. Define $\bS^{d-1}$ to be the unit sphere in $(\R^d,\|\cdot\|)$. For any measurable set $A \subset \bS^{d-1}$ and $0 < b < \delta$,  we have
\begin{align*}
    \bP(\theta_1 \in A \mid V_1 \leq b) &= \frac{\bP(\theta_1 \in A, V_1 \leq b)}{\bP(V_1 \leq b)} \\
    & \leq \frac{\bP(\theta_1 \in A, R_1 \leq (\frac{b}{c_d(1 - \varepsilon)f_0(x)})^{1/d} )}{\bP(R_1 \leq (\frac{b}{c_d(1 + \varepsilon)f_0(x)})^{1/d} )}\\ 
    & = \frac{\int_{\{(z - x)/\Vert z - x \Vert \in A,\Vert z - x \Vert \leq (\frac{b}{c_d(1 - \varepsilon)f_0(x)})^{1/d}\}}f_1(z) \dif z}{ \int_{\{\Vert z - x \Vert \leq (\frac{b}{c_d(1 + \varepsilon)f_0(x)})^{1/d}\}} f_1(z) \dif z} \\
    & \leq \frac{(1 + \varepsilon)f_1(x)}{(1 - \varepsilon)f_1(x)} \frac{\lambda(\{(z - x)/\Vert z - x \Vert \in A,\Vert z - x \Vert \leq (\frac{b}{c_d(1 - \varepsilon)f_0(x)})^{1/d}\})}{ \lambda \{\Vert z - x \Vert \leq (\frac{b}{c_d(1 + \varepsilon)f_0(x)})^{1/d}\}}\\
    & = \frac{(1 + \varepsilon)f_1(x)}{(1 - \varepsilon)f_1(x)} \frac{\mu_{d}(A)b/c_d(1 - \varepsilon)f_0(x)}{b/c_d(1 + \varepsilon)f_0(x) }\\
    & = \frac{(1 + \varepsilon)^2}{(1 - \varepsilon)^2}\mu_{d}(A).
\end{align*}
Here $\mu_d$ is the normalized Hausdorff measure on sphere $\bS^{d - 1}$ such that $\mu_d({\bS^{d - 1}}) = 1$.
Similarly, we have
\[
\bP(\theta_1 \in A \mid V_1 \leq b) \geq \frac{(1 - \varepsilon)^2}{(1 + \varepsilon)^2}\mu_{d}(A).
\]
Therefore we obtain
\[
\lim_{b \to 0+}\bP(\theta_1 \in A \mid V_1 \leq b) = \mu_d(A).
\]
By L'Hopital's rule, we have
\[
\lim_{b \to 0+}\frac{\bP(\theta_1 \in A, V_1 \leq b)}{\bP(V_1 \leq b)} = \lim_{b \to 0+}\frac{\frac{\partial}{\partial b}\bP(\theta_1 \in A, V_1 \leq b)}{\frac{\dif}{\dif b}\bP(V_1 \leq b)} = \lim_{b \to 0+}\bP( \theta_1 \in A \mid V_1 = b) = \mu_d(A).
\]
Since the above holds for any measurable set $A \subset \bS^{d - 1}$, then as $v \to 0+$, we have
\begin{align}\label{eq:6.3}
\theta_1 \mid V_1 = vv_1 \xrightarrow{\enskip d \enskip} \mathrm{Unif}(\bS^{d-1}).
\end{align}
Let $\tilde{\theta}_1$ and $\tilde{\theta}_2$ be two independent copies of $\mathrm{Unif}(\bS^{d-1})$. Then combining with \eqref{eq:6.1} we have
\[
\begin{aligned}
    & \quad \limsup_{v \to 0+}\bP(V \leq vt \mid V_1 = vv_1, V_2 = vv_2) \\
    & \leq \limsup_{v \to 0+}\bP\Big(c_d^{-1}S(v_1^{1/d}, v_2^{1/d}, \theta_1, \theta_2) \leq \frac{1 + \varepsilon}{1 - \varepsilon}t \mid V_1 = vv_1, V_2 = vv_2\Big) \\
    & \leq \bP\Big(c_d^{-1}S(v_1^{1/d}, v_2^{1/d}, \tilde{\theta}_1, \tilde{\theta}_2) \leq \frac{1 + \varepsilon}{1 - \varepsilon}t\Big).
\end{aligned}
\]
Similarly, combining with \eqref{eq:6.2} we have
\[
\liminf_{v \to 0+}\bP\Big(V \leq vt \mid V_1 = vv_1, V_2 = vv_2\Big) \geq \bP\Big(c_d^{-1}S(v_1^{1/d}, v_2^{1/d}, \tilde{\theta}_1, \tilde{\theta}_2 )\leq \frac{1 - \varepsilon}{1 + \varepsilon}t\Big).
\]
Since the above holds for arbitrary $\varepsilon > 0$, we then obtain
\[
\lim_{v \to 0+}\bP\Big(\frac{V}{v}\leq t \mid V_1 = vv_1, V_2 = vv_2\Big) = \bP\Big(c_d^{-1}S(v_1^{1/d}, v_2^{1/d}, \tilde{\theta}_1, \tilde{\theta}_2 )\leq t\Big),
\]
for all $t > 0$. Therefore as $v \to 0+$, the condition distribution converges:
\[
\frac{V}{v} \mid V_1 = vv_1, V_2 = vv_2 \xrightarrow{\enskip d \enskip} c_d^{-1}S(v_1^{1/d}, v_2^{1/d}, \tilde{\theta}_1, \tilde{\theta}_2 ).
\]
Denote the density of $c_d^{-1}S(v_1^{1/d}, v_2^{1/d}, \tilde{\theta}_1, \tilde{\theta}_2 )$ by $f_{S, v_1, v_2}(t)$. By linear transform of density functions, we have

\[
\lim_{v \to 0+}vf_{V \mid V_1, V_2}(v \mid vv_1, vv_2) = \lim_{v \to 0+}f_{\frac{V}{v} \mid V_1, V_2}(1 \mid vv_1, vv_2)  = f_{S, v_1, v_2}(1).
\]
For any sequence ${v_n}$ converges to 0, define 
\begin{align*}
&g_n(v_1, v_2) = P_{ijk}(1, v_1, v_2)(v_nf_{V \mid V_1, V_2}(v_n\mid v_nv_1, v_nv_2))\frac{f_{V_1, V_2}(v_nv_1, v_nv_2)}{f_V(v_n)/v_n}\\
{\rm and}~~~&g(v_1, v_2) = \frac{P_{ijk}(1, v_1, v_2) f_{S, v_1, v_2}(1)}{\alpha(d)}.
\end{align*}
Notice that $Q_{ijk}(v) \leq 1$ is uniformly bounded.  Therefore
\[
Q_{ijk}(v_n) = \int_{D(1)}g_n \leq 1.
\]
Since $g_n(v_1, v_2) \to g(v_1, v_2)$, invoking Fatou's Lemma yields

\[
\int_{D(1)}g \leq \liminf_{n \to \infty}\int_{D(1)}g_n \leq 1.
\]
Thus $g$ is integrable on $D(1)$. Define $D(1, \delta)$ to be
\[
D(1, \delta) = \Big\{(v_1, v_2) \in \bR^2: \delta \leq v_1, v_2 \leq 1 - \delta, v_1 + v_2 \geq 1 + \delta\Big\};
\]
specifically, $D(1) = D(1, 0)$. Then for any $\varepsilon > 0$, there exists some $\delta > 0$ such that
\begin{align}\label{eq:6.4}
\int_{D(1, \delta)} g \geq \int_{D(1)}g - \varepsilon ~~~{\rm and}~~~
\int_{D(1, \delta)} g_n \geq \int_{D(1)}g_n - \varepsilon.
\end{align}
Notice $g(v_1, v_2) < \infty$ for every $(v_1, v_2 ) \in D(1, \delta)$ and $D(1, \delta)$ is a compact subset on $\bR^2$.  By the continuity of $g(v_1, v_2)$, we have $g(v_1, v_2)$ is uniformly bounded on $D(1, \delta)$.  By dominated convergence theorem on this compact set, we then have
\[
\lim_{n \to \infty} \int_{D(1, \delta)}g_n  = \int_{D(1,\delta)}g.
\]
Combining \eqref{eq:6.4} we obtain
\[
\int_{D(1)}g \geq \lim_{n \to \infty}\int_{D(1)}g_n - 2\varepsilon.
\]
Thus we obtain
\[
\lim_{n \to \infty}\int_{D(1)}g_n = \int_{D(1)}g.
\]
This shows
\[
\begin{aligned}
\lim_{v \to 0+}Q_{ijk}(v) &= \lim_{v \to 0+}\int_{D(1)}P_{ijk}(1, v_1, v_2)(vf_{V \mid V_1, V_2}(v\mid vv_1, vv_2))\frac{f_{V_1, V_2}(vv_1, vv_2)}{(f_V(v)/v)} \dif v_1 \dif v_2 \\
&= \int_{D(1)}\frac{P_{ijk}(1, v_1, v_2) f_{S, v_1, v_2}(1)}{\alpha(d)} \dif v_1 \dif v_2 \\
&=: c_{ijk}(d).
\end{aligned}
\]
It can be seen that the constant $c_{ijk}(d)$ does not depend on the distribution of $\nu_1, \nu_0$.

~\\
\textbf{Step 2}: Next we consider the case of $d = 1$ and give the exact value of $c_{ijk}(1)$.

Following the same notation as in \textbf{Step 1}, define  
\[
R_1 = \Vert Z_1 - x\Vert,~~ \theta_1 =\frac{Z_1 - x}{\Vert Z_1 - x\Vert},~~ R_2 = \Vert Z_2 - x\Vert,~~ \theta_2 = \frac{Z_2 - x}{\Vert Z_2 - x \Vert}. 
\]
Since $d = 1$, then $Z_1, Z_2 \in \bR$ and $\theta_1, \theta_2 \in \{1, -1\}$. Therefore $V$ can be directly expressed by $V_1, V_2$ as follows
\[
V = \begin{cases}
       V_1 + V_2, & \mbox{ if } \theta_1 = -\theta_2,\\
    \max\{V_1, V_2\}, & \mbox{ if } \theta_1 = \theta_2.
    \end{cases}
\]
Obviously we have
\[
(V - V_1)(V - V_2)(V_1 + V_2 - V) = 0.
\]
Thus when $i, j, k > 0$, we have
\[
Q_{ijk}(v) = 0.
\]
Observe that
\[
\begin{aligned}
Q_{ijk}(v) & = \bE\Big[P_{ijk}(V, V_1, V_2) \mid V = v\Big] \\
            & = \bE\Big[P_{ijk}(V, V_1, V_2) \mid  \theta_1 = \theta_2, V  = v\Big]\cdot \bP(\theta_1 = \theta_2 \mid V = v)\\
            &~~~+ \bE\Big[P_{ijk}(V, V_1, V_2) \mid  \theta_1 = -\theta_2, V = v\Big]\cdot\bP(\theta_1 = -\theta_2 \mid V = v),
\end{aligned}
\]
and
\[
    \bP(\theta_1 = -\theta_2 \mid V \leq v) = \frac{\bP(\theta_1 = - \theta_2, V \leq v)}{\bP(V \leq v)}  = \frac{\bP(\theta_1 = - \theta_2, V_1 + V_2 \leq v)}{\bP(V \leq v)}.
\]
By \eqref{eq:6.3} in \textbf{Step 1}, we have $\theta_i \mid V_i = v$ converges in distribution to $\mathrm{Unif}(\bS^{d-1}) = \mathrm{Unif}\{1, - 1\}$. Letting $\tilde{\theta}_1, \tilde{\theta}_2$ be independent  $\mathrm{Unif}\{1, - 1\}$, we have
\[
\lim_{v \to 0+}\bP(\theta_1 = - \theta_2\mid V_1 + V_2 \leq v) = \bP(\tilde{\theta}_1 = -\tilde{\theta}_2) = 1/2.
\]
Accordingly, it holds true that
\[
    \bP(\theta_1 = -\theta_2 \mid V \leq v) = \frac{\bP(\theta_1 = - \theta_2, V_1 + V_2 \leq v)}{\bP(V \leq v)}  = \frac{\bP(V_1 + V_2 \leq v)}{2\bP(V \leq v)}.
\]
From Lemma \ref{lem5.1}, we have
\[
\lim_{v \to 0+}\frac{\bP(V \leq v)}{v^2} = \frac{\alpha(1)}{2}\Big(\dratio\Big)^2= \frac{3}{4}\Big(\dratio\Big)^2.
\]
We also have
\[
\lim_{v \to 0+} f_{V_1}(v) = \frac{f_1(x)}{f_0(x)} .
\]
Combining the above identities with the fact that $V_1, V_2$ are independent and have the same distribution, we obtain
\[
\begin{aligned}
    \lim_{v \to 0+}\frac{\bP(V_1 + V_2 \leq v)}{v^2} & = \lim_{v \to 0+}\int_{0 \leq x + y \leq v}\frac{f_{V_1}(x)f_{V_1}(y)}{v^2}\dif x\dif y  = \lim_{v \to 0+} \int_{0 \leq x + y \leq 1}f_{V_1}(vx)f_{V_1}(vy) \dif x\dif y \\
    & = \Big(\dratio\Big)^2\int_{0 \leq x + y \leq 1}\dif x\dif y = \frac{1}{2}\Big(\dratio\Big)^2.
\end{aligned}
\]
Therefore
\begin{align}\label{eq:6.5}
\lim_{v\to 0+} \bP(\theta_1 = -\theta_2 \mid V = v)  = \lim_{v\to 0+} \bP(\theta_1 = -\theta_2 \mid V \leq v) = \lim_{v \to 0+}\frac{\bP(V_1 + V_2 \leq v)/v^2}{2\bP(V \leq v)/v^2} = \frac{1}{3},
\end{align}
and similarly we obtain
\begin{align}\label{eq:6.6}
\lim_{v\to 0+} \bP(\theta_1 = \theta_2 \mid V = v) = \frac{2}{3}.
\end{align}

Recall that when $i, j, k > 0$, we have $Q_{ijk}(v) = 0$. We then consider the following two cases for $c_{ijk}(1)$.
~\\
\textbf{Case 1}: Only one of $i, j, k$ is $0$. Suppose $k = 0$ and $i,j > 0$. Then $\theta_1 = \theta_2$ implies $P_{ijk}(V, V_1, V_2) = 0$, and thus
    \[
    \begin{aligned}
    Q_{ijk}(v) &= \bE\Big[P_{ijk}(V, V_1, V_2) \mid  \theta_1 = -\theta_2, V = v\Big]\cdot \bP(\theta_1 = -\theta_2 \mid V = v) \\
    & = \bE\Big[P_{ijk}(V, V_1, V_2) \mid  \theta_1 = -\theta_2, V_1 + V_2 = v\Big]\cdot\bP(\theta_1 = -\theta_2 \mid  V = v) \\
    & = \bE\Big[V_1^j(v - V_1)^i/v^{i + j}\mid V_1 + V_2 = v, \theta_1 = -\theta_2\Big]\cdot\bP(\theta_1 = -\theta_2 \mid  V = v). \\
    \end{aligned}
    \]
    Since the density functions of $V_1, V_2$ converge to a constant around 0, thus the condistional distribution $(V_i/v) \mid V_i \leq v, \theta_i$ converges in disrtibution to $\mathrm{Unif}[0, 1]$. Letting $U_1, U_2$ be independent $\mathrm{Unif}[0, 1]$, then given $V_1 + V_2 = v, \theta_1 = - \theta_2$, we have
    \[
    \frac{V_1}{v} \mid V_1 + V_2 = v, \theta_1 = - \theta_2 \xrightarrow{\enskip d \enskip} U_1 \mid U_1 + U_2 = 1 \sim \mathrm{Unif}[0, 1], \quad \text{as}~~ v \to 0+.
    \]
    Therefore
    \[
    \lim_{v \to 0+} \bE\Big[V_1^j(1 - V_1)^i/v^{i + j}\mid V_1 + V_2 = v, \theta_1 = -\theta_2\Big] = \int_0^1x^j(1 - x)^i \dif x = \frac{i!j!}{(i + j + 1)!}.
    \]
    Combining the above with \eqref{eq:6.5}, we have
    \[
    c_{ijk}(1) = \lim_{v \to 0+}Q_{ijk}(v) = \frac{1}{3}\cdot\frac{i!j!}{(i + j + 1)!}.
    \] 
    
    Next suppose $j = 0$, $i, k > 0$. Then if either $\theta_1 = - \theta_2$ or $\theta_1 = \theta_2, V_1 > V_2$, it holds true that $P_{ijk}(V, V_1, V_2) = 0$, and thus
     \[
    \begin{aligned}
    Q_{ijk}(v) &= \bE\Big[P_{ijk}(V, V_1, V_2) \mid  \theta_1 = \theta_2,  V_1 \leq V_2, V = v\Big]\cdot \bP(\theta_1 = \theta_2, V_1 \leq V_2 \mid V = v) \\
    & = \bE\Big[P_{ijk}(V, V_1, V_2) \mid  \theta_1 = \theta_2, V_1 \leq V_2 = v\Big]\cdot \bP(\theta_1 = \theta_2, V_1 \leq V_2 \mid  V = v) \\
    & = \bE\Big[(v - V_1)^iV_1^k/v^{i + k}\mid \theta_1 = \theta_2, V_1 \leq V_2 = v\Big]\cdot \bP(\theta_1 = \theta_2, V_1 \leq V_2 \mid   V = v). \\
    \end{aligned}
    \]
    Similarly, we have
     \[
    \frac{V_1}{v} \mid V_1 \leq V_2 = v, \theta_1 =  \theta_2 \xrightarrow{\enskip d \enskip} U_1 \mid U_1 \leq U_2 = 1 \sim \mathrm{Unif}[0, 1], \quad \text{as} ~~v \to 0+.
    \]
    And by symmetry, we have
    \[
    \bP(\theta_1 = \theta_2, V_1 \leq V_2 \mid   V = v) = \frac{1}{2}\bP(\theta_1 = \theta_2 \mid   V = v).
    \]
    Then by \eqref{eq:6.6} we obtain
    \[
      c_{ijk}(1) = \lim_{v \to 0+}Q_{ijk}(v) = \frac{i!k!}{(i + k + 1)!}\cdot\frac{1}{2}\cdot\frac{2}{3} = \frac{1}{3}\frac{i!k!}{(i + k + 1)!}.
    \]
    Since $i, j$ are symmtric, then the case when $i = 0$, $j, k > 0$ is similar.
    
\textbf{Case 2}: Only one of $i, j, k$ is not $0$. Suppose $i = j = 0$, $k > 0$. Then $\theta_1 = -\theta_2$ implies $P_{ijk}(V,V_1, V_2) = 0$, and thus
    \begin{align*}
    Q_{ijk}(v) &= \bE\Big[P_{ijk}(V, V_1, V_2) \mid  \theta_1 = \theta_2, V = v\Big]\cdot\bP(\theta_1 = \theta_2 \mid V = v) \\
    & = \bE\Big[P_{ijk}(V, V_1, V_2) \mid  \theta_1 = \theta_2, \max\{V_1, V_2\} = v\Big]\cdot \bP(\theta_1 = \theta_2 \mid  V = v). 
    \end{align*}
    By symmetry, we have
    \[
    \bP(V_1 \leq V_2 = v\mid \theta_1 = \theta_2) = \bP(V_2 \leq V_1 = v\mid \theta_1 = \theta_2).
    \]
    It then implies
    \begin{align*}
      \quad \bE[P_{ijk}(V, V_1, V_2) \mid  \theta_1 = \theta_2, \max\{V_1, V_2\} = v] &= \frac{\bE[P_{ijk}(V, V_1, V_2) \ind \{\max\{V_1, V_2\} = v\} \mid  \theta_1 = \theta_2]}{\bP(\max\{V_1, V_2\} = v\mid \theta_1 = \theta_2)}\\
     & = \frac{2\bE[P_{ijk}(V, V_1, V_2) \ind \{V_1 \leq V_2 = v\} \mid  \theta_1 = \theta_2]}{2\bP(V_1 \leq V_2 = v\mid \theta_1 = \theta_2)} \\
     & = \bE[P_{ijk}(V, V_1, V_2) \mid  \theta_1 = \theta_2, V_1 \leq V_2, V_2 = v] \\
     & = \bE[V_1^k/v^k\mid  \theta_1 = \theta_2, V_1 \leq V_2 = v].
    \end{align*}
We then have
       \[
      c_{ijk}(1) = \lim_{v \to 0+}Q_{ijk}(v) = \frac{2}{3}\cdot\frac{1}{k + 1}.
    \]

    Suppose $j = k = 0$, $i > 0$. Then $\theta_1 = \theta_2, V_1 > V_2$ implies $P_{ijk}(V, V_1, V_2) = 0$. Thus
    \[
\begin{aligned}
Q_{ijk}(v)          =& \bE\Big[P_{ijk}(V, V_1, V_2) \mid  \theta_1 = -\theta_2, V  = v\Big]\cdot \bP(\theta_1 = -\theta_2 \mid V = v)\\
            &+ \bE\Big[P_{ijk}(V, V_1, V_2) \mid  \theta_1 = \theta_2, V_1 \leq V_2 = v\Big]\cdot\bP(\theta_1 = \theta_2, V_1 \leq V_2 \mid V = v) \\
             =& \bE\Big[(v - V_1)^i/v^i \mid \theta_1 = -\theta_2, V_1 + V_2 = v\Big]\cdot \bP(\theta_1 = -\theta_2 \mid V = v) \\
            &+ \bE\Big[(v - V_1)^i/v^i \mid \theta_1 = \theta_2, V_1 \leq V_2 = v\Big]\cdot \bP(\theta_1 = \theta_2, V_1\leq V_2 \mid V = v).
\end{aligned}
\]
Combining the above with results in \textbf{Case 1}, we obtain
  \[
      c_{ijk}(1) = \lim_{v \to 0+}Q_{ijk}(v) = \frac{1}{3}\cdot\frac{1}{i + 1} +  \frac{1}{3}\cdot\frac{1}{i + 1} = \frac{2}{3}\cdot\frac{1}{i + 1}.
    \]
By symmetry, a similar result holds for $i = k = 0$, $j > 0$.

Combining all results in \textbf{Case 1} and \textbf{Case 2} completes the proof.

\subsection{Proof of Lemma \ref{lem5.3}}

\nb{By the proof of Theorem \ref{thm3.1}, for any $x \in \cX$ we have
\begin{align*}
\bE\Big[n_0 \nu_1(\cA_M(x)) \mid  X_i = x, W_i = 0\Big] =  \bE\Big[\frac{n_0}{n_1}n_1\nu_1(\cA_M(x)) \mid X_i = x,  W_i = 0\Big] \\
= \bE\Big[\frac{n_0}{n_1}K_M(i) \mid X_i = x, W_i = 0\Big].
\end{align*}
Therefore we have
\[
\bE\Big[\frac{n_0}{n_1}K_M(i) \mid W_i = 0\Big] \leq \bE\Big[\frac{n_0^2}{n_1^2} \mid X_i = x, W_i = 0\Big] \bE\Big[K_M(i)^2 \mid X_i =x,  W_i = 0\Big].
\]
By the proof of Lemma 3 in \citet[Page 264]{abadie2006large}, we have 
\[
\bE\Big[K_M(i)^q \mid X_i = x, \mathbf{W}, W_i = 0\Big] \leq \sum_{k = 1}^q S  _k \Big(\frac{n_1}{n_0}\Big)^k
\]
for some constant $\{S_k\}_{k \geq 1}$ depending only on index $k$.
Also, by Lemma S.3 in \cite{abadie2016matching}, we have
\[
\bE\Big[\Big(\frac{n}{n_1}\Big)^r\Big] \leq C_r,
\]
for some constant $C_r$ depending only on $r$. Thus we obtain 
\[
\bE\Big[n_0 \nu_1(\cA_M(x)) \mid  X_i = x, W_i = 0\Big] \leq \bE\Big[\frac{n_0^2}{n_1^2} \mid W_i = 0\Big] \bE\Big[S_1\frac{n_1}{n_0} + S_2 \frac{n_1^2}{n_0^2} \mid W_i = 0\Big],
\]
which is upper bounded by a constant depending only through $C_1, C_2, S_1, S_2$. Therefore 
\[
\bE\Big[n_0 \nu_1(\cA_M(x)) \mid  X_i = x, W_i = 0\Big] 
\]
is uniformly bounded for all $n$ and $x \in \cX$.
Similarly, we have 
\[
\bE\Big[n_0^2\nu_1(\cA_M(x))^2 \mid X_i = x, W_i = 0\Big] 
\]
is uniformly bounded in $n$ and $x \in \cX$.
}

\subsection{Proof of Lemma \ref{lem5.4}}

\nb{
Fix $\delta>0$. By Assumption \ref{assump4.3}(i), the support $\cX$ is compact, so that there exists a finite collection of points $\{z_1,\ldots,z_N\}\subseteq\cX$ such that
\[
\cX\subseteq\bigcup_{\ell=1}^{N}B(z_\ell,\delta/4),
\]
where $N=N(\delta)$ depends only on $\delta$ and the geometry of $\cX$. Define the event
\[
E_n:=\Big\{\text{for every }\ell\in\{1,\ldots,N\},~\#\big\{k:\,W_k=0,\,X_k\in B(z_\ell,\delta/4)\big\}\geq M\Big\}.
\]

\textbf{Step 1: On the event $E_n$, $\mathrm{diam}(\cA_M(X_i))\leq\delta$.} Fix any $z\in\cA_M(X_i)$. By definition of the catchment area, $X_i$ is one of the $M$ nearest neighbors of $z$ among the control covariates $\{X_k\}_{W_k=0}$. Since $z\in\cX$, there exists some $\ell\in\{1,\ldots,N\}$ with $z\in B(z_\ell,\delta/4)$. By the diameter of $B(z_\ell,\delta/4)$, every point in $B(z_\ell,\delta/4)$ lies within distance $\delta/2$ of $z$, i.e.,
\[
B(z_\ell,\delta/4)\subseteq B(z,\delta/2).
\]
On $E_n$, the ball $B(z_\ell,\delta/4)$ contains at least $M$ control units, so that the $M$-th nearest control neighbor of $z$ lies within distance $\delta/2$. Consequently,
\[
\|X_i-z\|\leq\delta/2.
\]
This holds for every $z\in\cA_M(X_i)$, and hence $\cA_M(X_i)\subseteq B(X_i,\delta/2)$ and
\[
\mathrm{diam}\big(\cA_M(X_i)\big)\leq\delta.
\]

\textbf{Step 2: Show that $\bP(E_n)\to 1$.} By Assumption \ref{assump4.3}(ii), the density $f_X$ of $X$ is bounded below by $c_{\ell}>0$ on $\cX$. Recall that the density of $X$ given $W=0$ is 
\[
f_0(x)=f_X(x)(1-e(x))/(1-p). 
\]
By Assumption \ref{assump4.2}(ii), $1-e(x)\geq\eta$, so that
\[
f_0(x)\geq\frac{c_{\ell}\,\eta}{1-p}=:c_0>0,\qquad\text{for all }x\in\cX.
\]
For each $\ell$, let 
\[
V_\ell:=\nu_0(B(z_\ell,\delta/4)\cap\cX) 
\]
be the $\nu_0$-mass of the $\ell$-th cover ball restricted to the support, where $\nu_0$ denotes the law of $X\given W=0$. Since $z_\ell\in\cX$ and $\cX$ has nonempty interior at $z_\ell$ (or, more generally, since $\cX$ is the closure of a regular set; alternatively, after a possibly inflated cover, every $z_\ell$ is an interior point), the Lebesgue volume of $B(z_\ell,\delta/4)\cap\cX$ is bounded below by some $v_\delta>0$ depending only on $\delta$ and $\cX$, and hence
\[
V_\ell\geq c_0\,v_\delta=:\beta_\delta>0,\qquad \ell=1,\ldots,N.
\]

Conditional on $\mathbf{W}=(W_1,\ldots,W_n)$, the control covariates $\{X_k\}_{W_k=0}$ are i.i.d.\ from $\nu_0$, so that the count 
\[
N_\ell:=\#\{k:W_k=0,X_k\in B(z_\ell,\delta/4)\} 
\]
is $\mathrm{Binomial}(n_0,V_\ell)$. By the multiplicative Chernoff bound, for any $\ell$,
\[
\bP(N_\ell<M\mid \mathbf{W})\leq\bP\Big(N_\ell\leq\tfrac{1}{2}n_0V_\ell\,\Big|\,\mathbf{W}\Big)\leq\exp\Big(-\tfrac{1}{8}n_0V_\ell\Big)\leq\exp\Big(-\tfrac{1}{8}n_0\beta_\delta\Big),
\]
provided that $n_0$ is large enough so that 
\[
\frac{1}{2}n_0V_\ell\geq M,
\]
which holds for all sufficiently large $n$ since $\beta_\delta>0$ and $n_0/n\to 1-p$ almost surely. Taking a union bound over $\ell\in\{1,\ldots,N\}$ then yields
\[
\bP\big(E_n^{\,c}\,\big|\,\mathbf{W}\big)\leq N\exp\Big(-\tfrac{1}{8}n_0\beta_\delta\Big).
\]
Since $n_0\to\infty$ almost surely, the righthand side converges to $0$ almost surely. The dominated convergence theorem (the bound is at most $1$) then yields
\[
\bP(E_n^{\,c})=\bE\big[\bP(E_n^{\,c}\mid\mathbf{W})\big]\to 0,
\]
i.e., $\bP(E_n)\to 1$.

\vspace{0.2cm}

Combining {\bf Steps 1 and 2}, we have
\[
\bP\big(\mathrm{diam}(\cA_M(X_i))>\delta\big)\leq \bP(E_n^{\,c})\to 0,
\]
and since $\delta>0$ was arbitrary, $\mathrm{diam}(\cA_M(X_i))\xrightarrow{\bP}0$.
}

{
\bibliographystyle{apalike}
\bibliography{matching}

@article{han2022azadkia,
	Author = {Han, Fang and Huang, Zhihan},
	Journal = {The Annals of Applied Probability},
	Number = {6},
	Pages = {5172--5210},
	Title = {Azadkia-{C}hatterjee's correlation coefficient adapts to manifold data},
	Volume = {34},
	Year = {2024}}

@article{demirkaya2024optimal,
	Author = {Demirkaya, Emre and Fan, Yingying and Gao, Lan and Lv, Jinchi and Vossler, Patrick and Wang, Jingbo},
	Journal = {Journal of the American Statistical Association},
	Number = {545},
	Pages = {297--307},
	Publisher = {Taylor \& Francis},
	Title = {Optimal nonparametric inference with two-scale distributional nearest neighbors},
	Volume = {119},
	Year = {2024}}

@article{li2024matching,
	Author = {Li, Xuqiao and Yan, Ying},
	Journal = {arXiv preprint arXiv:2407.08468},
	Title = {Matching-Based Policy Learning},
	Year = {2024}}

@article{huo2023adaptation,
	Author = {Huo, Yiyi and Fan, Yingying and Han, Fang},
	Journal = {arXiv preprint arXiv:2311.16486},
	Title = {On the adaptation of causal forests to manifold data},
	Year = {2023}}

@article{lu2023flexible,
	Author = {Lu, Sizhu and Ding, Peng},
	Journal = {arXiv preprint arXiv:2305.17643},
	Title = {Flexible sensitivity analysis for causal inference in observational studies subject to unmeasured confounding},
	Year = {2023}}

@book{boots1999spatial,
	Author = {Okabe, A and Boots, B and Sugihara, K and Chiu, S},
	Edition = {2nd},
	Journal = {Geographical information systems},
	Pages = {503--526},
	Publisher = {Wiley},
	Title = {Spatial Tessellations: Concepts and Applications of Voronoi Diagrams},
	Year = {2000}}

@article{lin2024consistency,
	Author = {Lin, Ziming and Han, Fang},
	Journal = {Biometrika},
	Number = {1},
	Pages = {asag005},
	Title = {On the consistency of bootstrap for matching estimators},
	Volume = {113},
	Year = {2026}}

@book{ding2024first,
	Author = {Ding, Peng},
	Publisher = {CRC Press},
	Title = {A First Course in Causal Inference},
	Year = {2024}}

@article{holzmann2024multivariate,
	Author = {Holzmann, Hajo and Meister, Alexander},
	Journal = {arXiv preprint arXiv:2407.08494},
	Title = {Multivariate root-n-consistent smoothing parameter free matching estimators and estimators of inverse density weighted expectations},
	Year = {2024}}

@article{ulloa2024propensity,
	Author = {Ulloa-P{\'e}rez, Ernesto and Carone, Marco and Luedtke, Alex},
	Journal = {arXiv preprint arXiv:2409.19230},
	Title = {Propensity score augmentation in matching-based estimation of causal effects},
	Year = {2024}}

@article{cattaneo2023rosenbaum,
	Author = {Cattaneo, Matias D and Han, Fang and Lin, Zhexiao},
	Journal = {Biometrika},
	Number = {1},
	Pages = {asae062},
	Title = {On {R}osenbaum's Rank-based Matching Estimator},
	Volume = {112},
	Year = {2025}}

@article{he2024propensity,
	Author = {He, Yihui and Han, Fang},
	Journal = {Biometrika},
	Pages = {asae026},
	Publisher = {Oxford University Press},
	Title = {On propensity score matching with a diverging number of matches},
	Year = {2024}}

@article{imbens2024causal,
	Author = {Imbens, Guido W},
	Journal = {Annual Review of Statistics and Its Application},
	Pages = {123--152},
	Publisher = {Annual Reviews},
	Title = {Causal inference in the social sciences},
	Volume = {11},
	Year = {2024}}

@article{abadie2011bias,
	Author = {Abadie, A. and Imbens, G. W.},
	Journal = {Journal of Business and Economic Statistics},
	Pages = {1--11},
	Title = {Bias-corrected matching estimators for average treatment effects},
	Volume = {29},
	Year = {2011}}

@article{abadie2012martingale,
	Author = {Abadie, Alberto and Imbens, Guido W},
	Journal = {Journal of the American Statistical Association},
	Number = {498},
	Pages = {833--843},
	Publisher = {Taylor \& Francis},
	Title = {A martingale representation for matching estimators},
	Volume = {107},
	Year = {2012}}

@article{abadie2008failure,
	Author = {Abadie, Alberto and Imbens, Guido W},
	Journal = {Econometrica},
	Number = {6},
	Pages = {1537--1557},
	Publisher = {Wiley Online Library},
	Title = {On the failure of the bootstrap for matching estimators},
	Volume = {76},
	Year = {2008}}

@book{MR3445317,
	Author = {Biau, G\'{e}rard and Devroye, Luc},
	Doi = {10.1007/978-3-319-25388-6},
	Isbn = {978-3-319-25386-2; 978-3-319-25388-6},
	Mrclass = {62G08 (60D05 60E15 62G07 62G20 62H30 68T05 68T10)},
	Mrnumber = {3445317},
	Mrreviewer = {Christian Rau},
	Pages = {ix+290},
	Publisher = {Springer},
	Title = {Lectures on the Nearest Neighbor Method},
	Url = {https://doi.org/10.1007/978-3-319-25388-6},
	Year = {2015}}

@article{MR682809,
	Author = {Bickel, Peter J. and Breiman, Leo},
	Fjournal = {The Annals of Probability},
	Issn = {0091-1798},
	Journal = {The Annals of Probability},
	Mrclass = {60F05 (62G10)},
	Mrnumber = {682809},
	Mrreviewer = {P. R\'{e}v\'{e}sz},
	Number = {1},
	Pages = {185--214},
	Title = {Sums of functions of nearest neighbor distances, moment bounds, limit theorems and a goodness of fit test},
	Url = {http://doi.org/10.1214/aop/1176993668},
	Volume = {11},
	Year = {1983}}

@book{MR1324786,
	Author = {Billingsley, Patrick},
	Edition = {3rd},
	Publisher = {John Wiley and Sons},
	Title = {Probability and Measure},
	Year = {1995}}

@article{MR532236,
	Author = {Friedman, Jerome H. and Rafsky, Lawrence C.},
	Doi = {10.1214/aos/1176344722},
	Fjournal = {The Annals of Statistics},
	Issn = {0090-5364},
	Journal = {The Annals of Statistics},
	Mrclass = {62G10 (62H15)},
	Mrnumber = {532236},
	Mrreviewer = {B. Raja Rao},
	Number = {4},
	Pages = {697--717},
	Title = {Multivariate generalizations of the {W}ald-{W}olfowitz and {S}mirnov two-sample tests},
	Url = {https://doi.org/10.1214/aos/1176344722},
	Volume = {7},
	Year = {1979}}

@article{hahn1998role,
	Author = {Hahn, Jinyong},
	Journal = {Econometrica},
	Number = {2},
	Pages = {315--331},
	Publisher = {JSTOR},
	Title = {On the role of the propensity score in efficient semiparametric estimation of average treatment effects},
	Volume = {66},
	Year = {1998}}

@article{lin2023estimation,
	Author = {Lin, Zhexiao and Ding, Peng and Han, Fang},
	Journal = {Econometrica},
	Number = {6},
	Pages = {2187--2217},
	Title = {Estimation based on nearest neighbor matching: from density ratio to average treatment effect},
	Volume = {91},
	Year = {2023}}

@article{lin2022regression,
	Author = {Lin, Zhexiao and Han, Fang},
	Journal = {Journal of Econometrics},
	Pages = {106080},
	Title = {On regression-adjusted imputation estimators of the average treatment effect},
	Volume = {251},
	Year = {2025}}

@article{neyman1923applications,
	Author = {Neyman, Jersey},
	Journal = {Roczniki Nauk Rolniczych},
	Number = {1},
	Pages = {1--51},
	Title = {Sur les applications de la th{\'e}orie des probabilit{\'e}s aux experiences agricoles: Essai des principes},
	Volume = {10},
	Year = {1923}}

@article{rosenbaum1983central,
	Author = {Rosenbaum, Paul R and Rubin, Donald B},
	Journal = {Biometrika},
	Number = {1},
	Pages = {41--55},
	Publisher = {Oxford University Press},
	Title = {The central role of the propensity score in observational studies for causal effects},
	Volume = {70},
	Year = {1983}}

@article{shi2021ac,
	Author = {Shi, Hongjian and Drton, Mathias and Han, Fang},
	Journal = {Bernoulli},
	Number = {2},
	Pages = {851--877},
	Title = {On {A}zadkia-{C}hatterjee's conditional dependence coefficient},
	Volume = {30},
	Year = {2024}}

@article{stuart2010matching,
	Author = {Stuart, Elizabeth A},
	Journal = {Statistical Science},
	Number = {1},
	Pages = {1--21},
	Publisher = {NIH Public Access},
	Title = {Matching methods for causal inference: A review and a look forward},
	Volume = {25},
	Year = {2010}}

@article{voronoi1908nouvelles,
	Author = {Voronoi, Georges},
	Journal = {Journal f{\"u}r die reine und angewandte Mathematik (Crelles Journal)},
	Number = {134},
	Pages = {198--287},
	Publisher = {De Gruyter},
	Title = {Nouvelles applications des param{\`e}tres continus {\`a} la th{\'e}orie des formes quadratiques. Deuxi{\`e}me m{\'e}moire. Recherches sur les parall{\'e}llo{\`e}dres primitifs.},
	Volume = {1908},
	Year = {1908}}

@article{devroye2017measure,
	Author = {Devroye, Luc and Gy{\"o}rfi, L{\'a}szl{\'o} and Lugosi, G{\'a}bor and Walk, Harro},
	Journal = {Journal of Applied Probability},
	Number = {2},
	Pages = {394--408},
	Publisher = {Cambridge University Press},
	Title = {On the measure of {V}oronoi cells},
	Volume = {54},
	Year = {2017}}

@article{abadie2006large,
	Author = {Abadie, Alberto and Imbens, Guido W},
	Journal = {Econometrica},
	Number = {1},
	Pages = {235--267},
	Publisher = {Wiley Online Library},
	Title = {Large sample properties of matching estimators for average treatment effects},
	Volume = {74},
	Year = {2006}}

@article{rubin1974estimating,
	Author = {Rubin, Donald B},
	Journal = {Journal of Educational Psychology},
	Number = {5},
	Pages = {688--701},
	Publisher = {American Psychological Association},
	Title = {Estimating causal effects of treatments in randomized and nonrandomized studies.},
	Volume = {66},
	Year = {1974}}

@techreport{abadie2002simple,
	Author = {Abadie, A and Imbens, Guido},
	Institution = {National Bureau of Economic Research},
	Title = {Simple and Bias-Corrected Matching Estimators for Average Treatment Effects},
	Year = {2002}}

@article{abadie2016matching,
	Author = {Abadie, Alberto and Imbens, Guido W},
	Journal = {Econometrica},
	Number = {2},
	Pages = {781--807},
	Publisher = {Wiley Online Library},
	Title = {Matching on the estimated propensity score},
	Volume = {84},
	Year = {2016}}
}

\end{document}